\documentclass[dvips,ps]{imsart}
\RequirePackage[OT1]{fontenc} \RequirePackage{amsthm,amsmath}
\RequirePackage[colorlinks]{hyperref} \RequirePackage{hypernat}
\usepackage{graphicx}

\pubyear{2006} \volume{0} \issue{0} \firstpage{1} \lastpage{23}
\startlocaldefs \numberwithin{equation}{section}
\theoremstyle{plain}
\newtheorem{Theorem}{Theorem}[section]
\newtheorem{Corollary}{Corollary}[section]
\newtheorem{Lemma}{Lemma}[section]

\newtheorem{Remark}{Remark}[section]
\endlocaldefs\newdimen\AAdi%
\newbox\AAbo%
\font\AAFf=cmex10
\def\AArm{\fam0 \mathrm}%
\def\AAk#1#2{\setbox\AAbo=\hbox{#2}\AAdi=\wd\AAbo\kern#1\AAdi{}}%
\def\AAr#1#2#3{\setbox\AAbo=\hbox{#2}\AAdi=\ht\AAbo\raise#1\AAdi\hbox{#3}}%
\def\BBc{{\AArm C\AAk{-1.02}{C}\AAr{.9}{I}{\AAFf\char"3F}}}%
\def\BBe{{\AArm I\!E}}%
\def\BBn{{\AArm I\!N}}%
\def\BBp{{\AArm I\!P}}%
\def\BBr{{\AArm I\!R}}%
\def\BBz{{\AArm Z\!\!Z}}%
\def \indjk{j(p+1)+(k+1)}
\def \I{\hbox{\rm 1\hskip -3pt I}}

\def \sigs{(\sigma,\bar s(\xi))}

\begin{document}
\begin{frontmatter}
\title{Parametric estimation in  noisy blind deconvolution
model: a new estimation procedure} \runtitle{Parametric estimation
in noisy blind discrete deconvolution}

\begin{aug}
\author{\fnms{Emmanuelle} \snm{Gautherat}\ead[label=e1]{gauthera@ensae.fr}}
\and
\author{\fnms{Ghislaine} \snm{Gayraud}\ead[label=e2]{gayraud@ensae.fr}}

\address{E. Gautherat, CREST, Timbre J340, 3 av. P. Larousse, 92241 Malakoff Cedex and,
Economic Faculty de Reims, France, \printead{e1}\\
G. Gayraud, CREST, Timbre J340, 3 av. P. Larousse, 92241 Malakoff
Cedex and, LMRS, University de Rouen, France, \printead{e2}}

\runauthor{E. Gautherat et al.}


\end{aug}

\begin{abstract}
In a parametric framework, the paper is devoted to the study of a
new estimation procedure for the inverse filter and the level noise
in a complex noisy blind discrete deconvolution  model. Our
estimation method is a consequence of the sharp exploitation of the
specifical properties of the Hankel forms. The distribution of the
input signal is also estimated. The strong consistency and the
asymptotic distribution of all estimates are established. A
consistent simulation study is added in order to demonstrate
empirically the computational performance of our estimation
procedures.
\end{abstract}


\end{frontmatter}

\section{Introduction}
Let   $(Y_t)_{t \in \BBz}$ be the output process of an unknown
deterministic linear  time-invariant sequence $(u_t)_{t \in \BBz}$,
which is driven by an unobservable input sequence $(X_t)_{t \in
\BBz}$ added with a noise sequence $(\sigma_0\; W(t))_{t \in \BBz}$,
where $\sigma_0$ is an unknown level noise. In other words,
$(Y_t)_{t \in \BBz}$ is issued of  the noisy blind deconvolution
model defined by
\begin{eqnarray} Y_t & =&  (u \star X)_t + \sigma_0 W_t
  = \sum_{k \in \BBz} u_k X_{t-k} + \sigma_0 W_t, \; \forall t \in \BBz
\label{model}
\end{eqnarray}
where $(X_t)_{t\in \BBz}$ is assumed to be a complex discrete finite
valued input process and  the real-valued filter $(u_t)_{t \in \BBz}
\in l_1 (\BBz)$ is supposed to invertible.

Since the sequence
 $(u_t)_{t \in \BBz}$ is invertible,  it exists $\theta=(\theta_t)_{t \in \BBz}$ the
inverse filter of $(u_t)_{t \in \BBz}$, such that $\sum_{t \in \BBz}
u_t \theta_{k-t}=\I_{\{k=0\}}$. Note that if $(u_t)_{t \in \BBz}$
has a finite length,  the system is a noisy moving average and if
$(\theta_t)_{t \in \BBz}$ has a finite length,  the system is a
noisy autoregressive.
From $n$ observations $(Y_t)_{\{1\leq t \leq n\}}$, the  objective
is to restitute the distribution of the input process $(X_t)_{t\in
\BBz}$ which requires the estimation of the level of noise
$\sigma_0$ and the filter $(u_t)_{t \in \BBz}$. In digital signal
framework, that is, when $(X_t)_{t\in \BBz}$ is  a discrete valued
input,  the Bayesian theory combined with the MCMC methods is often
used
 to obtain the  a posteriori distribution of the
signal process $(X_t)_{t\in \BBz}$ (Liu \& Chen [1995];
Li \& Shedden [2001]).
 Here, we adopt the approach which consists in estimating the
 inverse filter $(\theta)_{t \in \BBz}$ instead of estimating the filter itself $(u)_{t \in \BBz}$.
 The problem of estimating
$(\theta_t)_{t \in \BBz}$ as well as the distribution of the input
process from the observed output data $(Y_t)_{t=1,\cdots,n}$ is
known as an identification problem. This identification problem when
the distribution of the input signal is discrete with a finite
number of possible values, is discussed in a number of papers:
 Li [1992, 1993, 1995,
1999, 2003], Li \& Mbarek [1997], Gamboa  \& Gassiat [1996] or
Gassiat \& Gautherat [1998, 1999] and Gautherat [1997, 2002]. In the
real and complex cases without noise, Li [1995] proposed an
estimation method for the inverse filter when the support points of
the input signal are supposed to be known. Gamboa \& Gassiat [1996]
in the non-noisy real case under a general setting (unknown
distribution of the input signal) proved the consistency of the
inverse filter estimate and that of the cardinality of the support's
points. Gassiat \& Gautherat [1998] extended in some sense the
previous paper in considering data with an additive noise and
proposed consistent estimates for  the support's points and also for
the level noise. Gassiat \& Gautherat [1999] studied the rate of
convergence of the signal estimate and the inverse filter estimate
in the parametric framework whereas Gautherat [2002] also
established asymptotic distribution for the points and their
corresponding mass of the signal distribution. To construct the cost
function, some of these authors explicitly used the alphabet of the
signal values (Li [2003, 1992, 1993, 1995, 1999]; Li \& Mbarek
[1997]), whereas Gamboa \& Gassiat [1996], Gassiat \& Gautherat
[1998, 1999] and Gautherat [1997, 2002] used only the cardinality of
the alphabet (non communicative situation).

Chen \& Li [1995], Gamboa \& Gassiat [1996, 1997a, 1997b], and
Gunther \& Swindlehurst [2000] showed that incorporating the finite
alphabet information into blind deconvolution procedures can greatly
improve the accuracy of the filter estimation and that of the signal
distribution in the non-noisy situation. Due to the judicious
utilization of the finite alphabet information, these methods
enjoyed a number of desirable properties in the non-noisy situation
such that, the ability to handle (with super statistical
efficiently) a large class of filters $(u_t)_{t \in \BBz}$, the
ability to handle non-stationary (see the papers of Li) or non
i.i.d. signals (Gamboa  \& Gassiat [1996]) without modeling or
estimating their statistical characteristics.

In the model (\ref{model}), the blind deconvolution of the data is
much more complicated due to the presence of the noise. Gassiat and
Gautherat [1998, 1999] and Gautherat [1997, 2002] proposed  a
consistent estimation procedure which is based on the minimization
of a penalized  empirical constrast function. In practice, this
method needs  to adjust the penalty term and requires a starting
point which is near enough to the true value in order to avoid a
local minimizer.

Among these works, the main contribution of our paper is to provide
in the complex case, a new estimation procedure of the level of noise and
the inverse filter, which is a consequence of the sharp study of the
Hankel matrix when the noise is supposed to be Gaussian.
Our estimation procedure is built on an explicit empirical criteria (it is not the same in the papers
of Gassiat $\&$ Gautherat [1998, 1999] and Gautherat [1997, 2002]),  and it is based on the roots of an empirical function.
Our estimation method is competitive from a theoretical
point of view (consistency and asymptotic distribution of all estimates) and from a practical side (our numerical results are quite good).

The paper is organized as follows. The assumptions on the model
defined by (\ref{model}) are given in Section \ref{Charac:Sec}. In
the same section, the level of noise, the inverse filter and the law
of the input signal are characterized. These characterizations are used
to define in Section \ref{Esti:Sec} our estimation procedures. The
strong consistency and the asymptotic distribution of all the
estimates are stated in Section \ref{Resu:Sec}. The proofs are
postponed in Section \ref{Proof:Sec}. Section \ref{Simu:Sec} deals
with a simulation study in which the computational performance of
our estimation procedures is empirically demonstrated. Some
concluding remarks including a comparison of our numerical results with those given in Gassiat $\&$ Gautherat [1998], are made in Section \ref{Conclusion:Sec}.
\section{Assumptions and Characterization}\label{Charac:Sec}
%
%
\subsection{Assumptions}
\begin{itemize}
\item[{\bf (M1)}]
$(X_{k})_{k\in \BBz}$ is a sequence of discrete complex random
variables with a common discrete support $a=(a_{1},\ldots,a_{p}) \in
\BBc^p$; $a$ is unknown, $a_i\neq a_j$ for $i\neq j$.  The integer $p \geq 2$
is known. The components $(a_j)_{1\leq j \leq p}$ of $a$ are given
by the lexicographical order.

\item[{\bf (M1b)}]
$(X_k)_{k \in \BBz}$ is identically distributed with support points
$a \in \BBc^p$ and $\Pi=(\pi_j)_{\{j=1,\ldots,p \}}$ is such that
$\;\BBp (X_k= a_j) = \pi_j \; \forall j \in \{1,\ldots,p\} $ and for
any $k \in \BBz$. The integer $p \geq 2$ is known. The components
$(a_j)_{1\leq j \leq p}$ of $a$ are given by the lexicographical
order.
\item[{\bf (M2)}] $(X_{k})_{k\in \BBz}$ is a stationary
ergodic process.
\item[{\bf (M3)}]
$ \forall n \in \BBn^{\star}$, $\forall (j_{1},\ldots,j_{n}) \in
\{1,\ldots,p\}^n$, $\BBp
(X_{1}=a_{j_{1}},\ldots,X_{n}=a_{j_{n}})\;>\;0.$
\item[{\bf (M4)}] $\forall k \in \BBz$, $W_k=W^R_k + i W^I_k$, where $W^R=(W^R_k)_{k \in \BBz}$ and
$W^I=(W^I_k)_{k \in \BBz}$ are independent. The random variables
$(W^R_k)_{k\in \BBz}$ and $(W^I_{k})_{k\in \BBz}$ are real centered
i.i.d. Gaussian with a variance equal to $(1/2)$ and they are both
independent of $(X_{k})$.
\item[{\bf (M5)}] $U(x)=\sum_{k} u_{k} e^{ikx}$ is
continuous  and does not vanish on $[0,2\pi]$.
\end{itemize}
Note that Assumption {\bf (M5)} guarantees that $(u_t)_{t \in \BBz}$
is invertible and that  both $(u_t)_{t \in \BBz}$ and its inverse
are in $l_1(\BBz)$. Since $(u_t)_{t \in \BBz}$ is invertible, the
initial observed process $(Y_t)_{t \in \BBz}$ can be transformed by
applying any filter $s$ to $Y_t$. The resulting process $Z_t(s)$ is
then defined by
\begin{eqnarray}
\forall t \in \BBz,  Z_t(s)=(s \star Y)_t. \label{ConvZ}
\end{eqnarray}
\begin{itemize}
\item[{\textbf{(M6)}}]  The set of the filters is defined by $\Theta=\{ s(\xi) \in l_1(\BBz),
\xi \in {\cal K} \}$, where the function $s$ is known, $s \in {\cal
 C}^1$ is injective. The true inverse filter $\theta$ is in $\Theta$. The set ${\cal K}  \subset \BBr^d, \; d \in
\BBn^*$ is compact  such that if $\xi,{\tilde \xi} \in {\cal K}$
satisfy $s_k(\xi)=r s_{k-l}({\tilde \xi}), \; \forall k \in \BBz
\Rightarrow r=1, \; l=0, \mbox{ {\rm and} } \xi = {\tilde \xi}$. The
parameter $\xi$ is unknown.
\end{itemize}
\begin{Remark}
$\;$

\noindent (i) Assumption \textbf{(M6)} describes the set $\Theta$ at
which the filters we consider, belong.  The set $\Theta$ is
parametric via the unknown parameter vector $\xi$, so that
estimating the inverse filter $\theta$ is reduced to estimate
$\xi_0$.

\noindent (ii) Since the true inverse filter $\theta$ belongs to
$\Theta$ and  since the function $s$ is injective, it implies that
it exists $\xi_0$ in the interior of ${\cal K}$ such that
$s(\xi_0)=\theta$. Moreover, the last part of Assumption
\textbf{(M6)} guarantees the identifiability of the model; in
particular, it allows to avoid problems of scale and delay.

\noindent (iii) In the non noisy case, note that $Z_t(\theta)=X_t,
\; \; \forall t \in \BBz$ where $Z_t(\cdot)$ is defined by equation
(\ref{ConvZ}). Our estimation procedures will be based on the process
$Z_t(s)$, $s \in \Theta$.
\end{Remark}

\subsection{Characterizations} \label{charac}
In the non-noisy framework ($\sigma_0=0$), Gamboa \& Gassiat [1996]
stated that, under {\bf (M1)}, {\bf (M3)}, {\bf (M4)} and {\bf
(M5)}, $Z_t(s)$ takes at most $p$ distinct values if and only if
$s=\theta$ up to scale and delay. The definition of $\Theta$ as a
subset of $l_1(\BBz)$ avoids identifiability problems which could be
generated by possible changes in scale or delay. But since the
characterizations of $\sigma_0$ and $\theta$ are valid for any
inverse filter in $l_1(\BBz)$,  from now and only for this
paragraph, let us consider any filter $s$ in $l_1(\BBz)$. \\

The random variable $Z_t(s)$ for $s \neq \theta$ having at least $p$
points of support,  the characterization of $\sigma_0$ and $\Theta$
is made via a contrast function which is able to distinguish
discrete random variables whose support is of cardinality $p$ from
others whose support is of cardinality greater than $p$.\\
 Let
$d(s)$  be the conjugate moment column vector of dimension $(p+1)^2$
defined
 by
\begin{eqnarray}
 d_{\indjk}(s) = \BBe  \left(  (Z_0(s))^k \overline{(Z_0(s))^j}\right), \forall (j,k)\in \{0,\ldots,p\}^2, \forall s \in l_1(\BBz).\label{momentConj}
\end{eqnarray}
Note that
if one rewrites $ d(s)$ in a
$\left((p+1) \times (p+1)\right)$--matrix $ D(s)$ where $j$ and
$k$ correspond respectively to the columns and the rows then $ D(s)$ is a Hankel matrix.
For all $s \in l_1(\BBz)$ and $\sigma \geq 0$, it is always possible
to derive the conjugate pseudo-moment $\tilde d(\sigma,s)$
 from $d(s)$ in inverting the following
system:
\begin{eqnarray}
 d(s) = A (\sigma
\|s\|_2){\tilde d}(\sigma, s), \label{EquA}
\end{eqnarray}%
where $\|s\|_2$ denotes the $l_2$-norm of $s$. The matrix $A(\sigma
\|s\|_2)$, which depends on $\sigma \|s\|_2$, is an invertible
$(p+1)^2 \times (p+1)^2$-matrix such that for any pair $(j,k)\in
\{0,\ldots,p\}\times \{0,\ldots,p\}$ and $\forall \;  0 \leq m \leq
j, \;0 \leq l \leq k$, $ \forall \beta \in \BBr,$
$$A_{j(p+1)+k+1, m(p+1)+l+1}(\beta)=C_k^lC_j^m \gamma_{j-m,k-l}
\beta^{k-l+j-m}, $$
where $\gamma_{m,l}= \BBe ((W_0)^l  (\overline{W_0})^m)$. Note that
under {\bf(M4)} one gets $ \gamma_{m,l} = \I_{m=l} m!$. Since
$\gamma_{m,l}\neq 0 \Longleftrightarrow l=m$, denote
$\gamma_{m,l}=\gamma_m$. The inversion of the system (\ref{EquA}) is
explicit and $A^{-1}(\sigma \|s\|_2)$ is defined by the following:
$\forall \; (j,k)\in \{0,\ldots,p\}^2 \mbox{ {\rm such that } } \; 0\leq m\leq (j-k\vee 0),$ \mbox{ {\rm and } }
$\forall \beta \in \BBr_+$,
\begin{equation} \label{EquAm}A^{-1}_{j(p+1)+k+1,
m(p+1)+k-j+m+1}(\beta)=(-1)^{j-m}C_k^{j-m}C_j^m \gamma_{j-m}
\beta^{2(j-m)}.
\end{equation}
\begin{Remark} \label{pseudomoment} Denote by
$(R(s))_t= (s \star u \star X)_t$ and  $(V(s))_t= (s \star W)_t $,
$\forall t \in \BBz$ and $\forall s \in l_1(\BBz)$. Then,
 $\forall \; \sigma \leq \sigma_0$ and due to the
infinite divisibility of the Gaussian distribution, $\forall (j,k)
\in \{0,\ldots,p\}^2, \; {\tilde d}_{j(p+1) +(k+1)} (\sigma, s)$ is
equal to
$$  {\small \BBe\left[ \left((R(s))_0 +
\sqrt{\sigma_0^2 - \sigma^2 }(V(s))_0)\right)^k
\left((\overline{R(s)})_0 + \sqrt{\sigma_0^2 - \sigma^2 }
\overline{(V(s))}_0\right)^j\right].}$$
On the other hand, if $\sigma > \sigma_0$, $\tilde d(\sigma,s)$ has no
explicit form, one does not know if it corresponds to a random
variable moment. This explains why ${\tilde d}(\sigma,s)$ is called the
pseudo-moment vector.
\end{Remark}
Next, transform the pseudo-moment vector ${\tilde d}(\sigma,s)$ in a
$\left((p+1) \times (p+1)\right)$--matrix $\tilde D(\sigma,s)$ where $j$ and
$k$ correspond respectively to the columns and the rows i.e.
${\tilde D}(\sigma,s)= \left( {\tilde d}_{\indjk}(\sigma, s)
\right)_{(j,k) \in \{0,\ldots,p\}^2}$.  Then, let $J$ be the real
function defined by:
\begin{eqnarray}
\forall \sigma \geq 0,\; \; \forall s \in l_1(\BBz),\;\;\;J
(\sigma,s) &=& \det ({\tilde D}(\sigma,s)). \label{optimi}
\end{eqnarray}
The  characterizations of   $\theta$, $\sigma_0$, $a$ and $\Pi$ in
the model described  by relation  (\ref{model}) or equivalently by
relation (\ref{ConvZ}), are made through the function $J$ defined by
(\ref{optimi}). They are established in Gassiat \& Gautherat [1999]
for $\theta$ and $\sigma_0$, and in Gautherat [1997, 2002] for $a$
and $\Pi$ under a more general setting:
\begin{itemize}
\item (i) Under assumptions {\bf (M1)}, {\bf (M3)}, {\bf (M4)} and {\bf
(M5)},  the true level noise $\sigma_0$ and the true inverse filter
$\theta=s(\xi_0)$ satisfy
\begin{eqnarray}
\forall \sigma < \sigma_0,  & & J (\sigma, s(\xi)) >0, \forall \xi
\in
{\cal K} \label{CondBORD}\\
J(\sigma_0,s(\xi))=0   & iff &  s(\xi)=\theta \mbox{ {\rm up to
scale an delay}}\label{caractthetasig}.
\end{eqnarray}

\item (ii) Under {\bf (M1)}, the distribution points
$a=(a_i)_{\{i=1,\ldots,p\}}$ are the roots of the polynomial
function $p_{v^{\star}}$ in $\BBc [X]$  defined by
$p_{v^{\star}}(x)=\displaystyle{\sum_{j=0}^p }v^{\star}_j x^j$,
where $v^{\star}$ denotes the eigenvector associated with the
smallest eigenvalue of the matrix ${\tilde D}(\sigma_0,\theta)$.
\item (iii) Under   {\bf (M1b)}, the distribution $\Pi=(\pi_1,\ldots,\pi_p)$ is uniquely
determined as the solution of the following linear system in
$(q_1,\ldots,q_p)\in [0,1]^p$: $\BBe (X^k_0)=\sum_{i=1}^p q_i
a_i^k \; \; \forall k=0,\ldots,p-1.$
\end{itemize}
Note that (i) implies that $\sigma_0 = \min\{ \sigma >0; \exists s
\in l_1(\BBz) : \; J(\sigma,s)=0\}$.

\section{Estimation procedures}\label{Esti:Sec}
To construct our estimates, we consider again the filters of the form
$s(\xi)$, where $\xi \in {\cal K}$ is unknown and $s$ is a known
function. For any $\xi \in {\cal K}$, let us consider the truncated
sequence $\bar s(\xi)$ of $s(\xi) \in \Theta \subset l_1(\BBz)$ as $
\bar s_k(\xi)= s_k(\xi) \I_{|k| \leq k(n)} \; \forall k \in \BBz$,
where $k(n)$ is a sequence of nonnegative integers increasing with
$n$. Denote also $\|s(\cdot)\|_{n,2}=\|\bar s(\cdot)\|_{2}$. Define
the truncated version of relation (\ref{ConvZ})
that is,
$$
\forall \xi \in {\cal K}, \; \; \forall t=1+k(n),\ldots,
n-k(n),\;\;\;Z_{t}( \bar s(\xi))=\sum_{k=-k(n)}^{k(n)} s_k(\xi)
Y_{t-k}.
$$
Denote $d_n(\xi)$ the empirical conjugate moment vector of
dimension $(p+1)^2$, whose general term $ d_{j(p+1)+k+1,n}(\xi)$
is defined as the empirical version of (\ref{momentConj}),
$$
\forall \xi \in {\cal K},\;\;\;d_{j(p+1)+k+1,n}(
\xi)=\frac{1}{n-2k(n)} \sum_{t=1+k(n)}^{n-k(n)} Z_{t}^k( \bar
s(\xi)) \overline{Z_{t}^j(\bar s(\xi))}.
$$
Then, similarly  to (\ref{EquA}), define $\tilde d_n( \xi)$  as
the empirical conjugate pseudo-moment vector  whose general term
$\tilde d_{j(p+1)+k+1,n}( \xi)$ is the solution of the following
triangular system:
\begin{eqnarray}
\forall \xi \in {\cal K},\forall \sigma \geq 0,\;\;\; \tilde
d_n(\sigma, \xi) &=&A^{-1}(\sigma \| \bar s(\xi)\|_{2} ) d_n(\xi),
\label{EmpiriquePseudoMoment}
\end{eqnarray}
where $A^{-1}(\beta)$ is the matrix defined by (\ref{EquAm}).
Finally, let $J_n$ be the empirical version of  $J$ defined by
(\ref{optimi}):
$$ \forall \xi \in  {\cal K},\forall \sigma
\geq 0, \; \; J_n (\sigma, \xi)= \det (\tilde D_n\sigs),$$
where $\tilde D_n$ is the $((p+1) \times (p+1))$-matrix
corresponding to the rewriting of the  empirical pseudo-moment
vector $\tilde d_n(\sigma, \xi)$ in a matrix form.
Finally, let us define all the estimates:
\begin{itemize}
\item
$( \widehat{\sigma}_{0},\widehat{\xi}_0 \;)$ is the  solution of the
following system
$$ \left\{\begin{array}{l}
 J_n(\widehat{\sigma}_{0},
\widehat{\xi}_0)=0 , \\
\widehat{\sigma}_{0}= \min \left\{ \sigma \in \BBr_+; \exists \xi
\in {\cal K}: \;  J_n(\sigma, \xi)=0\right\}.
\end{array} \right.$$
\item $\widehat{\theta}$ is
defined as $\widehat{\theta}=\bar s(\widehat{\xi}_0)$.
\item The support points  $(\widehat{a}_1,\ldots,\widehat{a}_p)$ in
$\BBc^p$ are the roots rearranged by the lexicographic order  of the
polynomial function
$p_{\widehat{v}^{\star}}(x)=\displaystyle{\sum_{j=0}^p}
\widehat{v}_j^{\star} x^j$ in $\BBc [X]$, where
$\widehat{v}^{\star}=(\widehat{v}_0^{\star},\ldots,\widehat{v}_p^{\star})$
denotes the eigenvector associated with the smallest eigenvalue of
the matrix ${\tilde D}_n (\widehat{\sigma}_{0},\bar
s(\widehat{\xi}_0))$.
\item  The probability vector $(\widehat{\pi}_1,\ldots,\widehat{\pi}_p)$ is uniquely
determined as the solution of the following linear system in
$(q_1,\ldots,q_p)\in [0,1]^p$:
$ \tilde d_{k+1,n} (\widehat{\sigma}_{0},\bar
s(\widehat{\xi}_0))=\sum_{i=1}^p q_i \widehat{a}_i^j, \; \forall
j=0,\ldots,p-1.$
 \end{itemize}
{\bf  Existence of
$(\widehat{\sigma}_{0},\widehat{\xi}_0\;)$ for any $n \in
\BBn^{*}$.} \\
For any $n \in \BBn^{\star}$, fix $\xi \in {\cal K}$ and observe
$d_{n}(\xi)$. Then,  ${\tilde d}_{n}(\sigma,\xi)$ is obtained via
equation (\ref{EmpiriquePseudoMoment}) and depends only on the
unknown parameter $\sigma$  throughout $\sigma \|\bar s(\xi)\|_{2}$
since ${\bar s}(\xi)$ is supposed to be fixed.
 For any $v \in \BBc^{p+1}_{\star}$, consider the
hermitian form
$$Q_v(\sigma)= \sum_{k=0}^{p} \sum_{j=0}^{p} {\tilde d}_{j(p+1) +
(k+1),n}(\sigma,\xi)  v_k \overline{v}_j.$$
Note that $Q_v$ is a polynomial function in $\sigma$ with real
coefficients. For any $v \in \BBc^{p+1}_{\star}$, note that
$Q_v(0)=Det(D_n(\xi)) >0$, where $D_n(\xi)$ is a Hankel matrix of dimension $((p+1)\times
(p+1))$ which corresponds to the rewriting of the vector
$d_n(\xi)$. Note also that the highest degree of $Q_v$ is  equal to
$2p$.
\\
Then, if $p$ is an odd number, the coefficient of the highest term
is equal to $(-1)^{p}\; p! (\|s(\xi)\|_{n,2})^{2p} v_p {\overline
v}_p$ which is negative. It entails that it exists a real
positive ${\tilde \sigma}_n$ such that $Q_v({\tilde \sigma}_n)=0$
that is, the hermitian form is degenerate and then
$J_n(\tilde{\sigma}_{n}, \xi)=0$. \\
On the other hand, if $p$ is an even number, choose any
$v=(v_0,v_1,\ldots,v_p)'\in \BBc^{p+1}_{\star}$ such that $v_p=0$,
where the $v'$ denotes the transposed vector of $v$. Then, the term
of the highest degree is negative and is of order $\sigma^{2(p-1)}$.
As previously, it allows to conclude to the existence of a real
positive ${\tilde \sigma}_n$ which is a zero for $J_n(\cdot, \xi)$.
\section{Main results}\label{Resu:Sec}
Up to now,  for any function  $F$ with one or more arguments, set
$D^{r}F(y)$ be the value at $y$ of the $r$-th differential of $F$
and set $\partial_{i_1,\cdots,i_l}^{r_{i_1},\cdots,r_{i_l}}F(y)$ be
the value at $y$ of the $(\displaystyle{\sum_{k=1}^l} r_{i_k})$-th
partial derivative of $F$, where $r_{i_k}$ is the order of the
derivative with respect to its $i_k$-th coordinate. Denote also
$d(\xi)$ the $(p+1)^2$-vector with components $d_{j(p+1) +
(k+1),n}(\xi):=d_{j(p+1) +
(k+1),n}(s(\xi))$
defined by relation (\ref{momentConj}). \\
Some extra other assumptions are needed to establish the consistency
and the asymptotic distribution of all estimates.
\begin{itemize}
\item[{\bf (M7)}]$k(n)=o({\sqrt n}),\;\;\;
\displaystyle{\sum_{|k|>k(n)}} |\theta_{k}|=o(\frac{1}{\sqrt{n}})
\mbox{ {\rm as } }n \; \rightarrow \; +\infty$
\item[{\bf (M8)}] $\sqrt{n}\left(
d_n(\xi_0)-d(\xi_0), D^{1}d_{n}(\xi_0)- D^{1}d(\xi_0)\right)
\xrightarrow[n \; \rightarrow \; + \infty]{{\cal L}}
 {\cal N}(0,\Gamma)$. Denote
 $\Gamma_1$ the asymptotic variance of $\sqrt{n}\left( d_{n}(\xi_0)-d(\xi_0)\right)$.
\item[{\bf (P)}] The application $\xi \; \mapsto  \; s(\xi)$ is twice continuously
differentiable. For any $i=1,\ldots, d$, $(\partial_i^1
s_{k}(\xi_0))_{k\in\BBz}$ and $(\partial_{i}^2
s_{k}(\xi_0))_{k\in\BBz}$ are in $l_{1}(\BBz)$.
 Moreover,
$\left((\partial_1^1 s_{k}(\xi_0))_{k\in\BBz},\ldots, (\partial_d^1
s_{k}(\xi_0))_{k\in\BBz}\right)$ and $(s_k(\xi_0))_{k\in\BBz}$ are
linearly independent.
\end{itemize}
\begin{Theorem}
\label{consistance} Suppose that assumptions {\bf (M1)}-{\bf (M7)}
hold, then as $n$ goes to infinity, $\widehat{\sigma}_{0}$ converges
a.s. to $\sigma_0$ and $\left\| {\widehat \xi}_0- \xi_0  \right\| $
converges a.s. to  0,  where $\| \cdot \|$ denotes the Euclidean
norm
 in $\BBr^d$.
\end{Theorem}
\begin{Corollary}
 \label{consistanceapi} Suppose that assumptions {\bf (M1b)}, {\bf (M2)}-{\bf
(M7)} hold, then as $n$ goes to infinity, both $\left\| {\widehat
a}- a \right\|$ and $\left\| {\widehat \Pi}- \Pi \right\|$ converge
a.s. to 0,  where $\| \cdot \|$ denotes the Euclidean norm
 in $\BBr^p$.
\end{Corollary}
Before giving the asymptotic distribution of our estimates, recall
that the vector $v^* \in \BBc^{p+1}$ is the eigenvector associated
with the smallest eigenvalue of ${\tilde D}(\sigma_0,\theta)$, and
denote $v(b)$ any vector in $\BBc^{p+1} \cap \{ \|\cdot \|_2=1\}$
associated with $b \in \BBc^p$ such that the $v(b)_j$'s are the
complex coefficients of the polynomial function in $\BBc [X]$ having
the components of $b$ as roots. In particular, note that $v(a)=v^*$.
\begin{Theorem}
\label{propspeed} Under assumptions  {\bf (M1)}-{\bf (M8)} and {\bf
(P)},
 \begin{eqnarray*}
 \sqrt{n}\left(\begin{array}{c} \widehat{\xi}_0 - \xi\\
 \widehat{\sigma}_{0}-\sigma_0 \end{array}\right) & \xrightarrow[n \; \rightarrow \; +\infty]{{\cal L}}&
 {\cal N}\left(0_{d+1},N D^1 h(\tilde d(\sigma_0,\xi_0))  M(D^1 h(\tilde
d(\sigma_0,\xi_0)))' N' \right),
 \end{eqnarray*}
where $0_{d+1}$  {\rm is the $(d+1)$-dimensional zero vector  and }
\begin{eqnarray*}
 M& =&  A^{-1}(\sigma_0 \|\theta\|_2) \Gamma_1
(A^{-1}(\sigma_0 \|\theta\|_2))' ,\\
N & =& \frac{1}{\alpha}\left(
\begin{array}{c}
 \partial_2^2 J(\sigma_0,\xi_0)^{-1} \;
\partial^{1,1}_{1,2} J(\sigma_0,\xi_0)
\\
1
\end{array}\right),\\
\alpha&=&- \partial_1^1 J(\sigma_0,\xi_0)\\
 h (\tilde d(\sigma_0,\xi_0))& =& det\left( \left( \tilde d_{j(p+1)+i+1}(\sigma_0,\xi_0)\right)_{i,j\in \{0,\cdots,p\}}\right).\\
\end{eqnarray*}
\end{Theorem}
\begin{Corollary}
\label{propspeedapi} Under {\bf (M1b)}, {\bf (M2)}-{\bf (M8)} and
{\bf (P)} and using the previous notations, one gets
\begin{eqnarray*}
\sqrt{n} ({\hat a} - a) & \xrightarrow[n \; \rightarrow \;
+\infty]{{\cal L}}&
 {\cal N}\left(
0_{p}, \frac{1}{4|v_p^{\star}|^4} C^{-1} \; B \; R \; M {\bar B}'\;
{\bar R}' C^{-1}
 \right), \\
\sqrt{n} ({\hat \Pi} - \Pi) & \xrightarrow[n \; \rightarrow \;
\infty]{{\cal L}}&
 {\cal N}\left(
0_{p}, G\;R \; M  \;{\bar R}'\; {\overline G'}
 \right),
 \end{eqnarray*}
where,
\begin{eqnarray*}
C&=& {\rm diag}(K_1,\ldots,K_p), \mbox{{\rm with }} K_j= \BBe \; (
\prod_{i=1, i \neq j}^p |X_0 -a_i|^2),\\
B && \mbox{{\rm is a $p\times (p+1)^2$-matrix which is defined by
its columns as follows }} \\
&&\forall l \in \{1,\ldots,p\},\; \mbox{{\rm with }}  (i,j) \in
\{0,\ldots,p\}^2,\\
 &&\;\; B_{l,i(p+1)+j+1} =\left( \overline{
\partial_l^1 v_i (a)} \; v_j^* +
\partial_l^1 v_j (a){\overline v_i^*}\right) \\
R& =&  \left( {\rm Id}_{(p+1)^2} + \partial_{1,2}^{1,1}\tilde
d(\sigma_0,\xi_0) D^1 h( \tilde d(\sigma_0,\xi_0))  N
\right) \\
&& \mbox{\rm where $\mbox{{\rm Id}}_d$ is the identity matrix of size $d$,}\\
G& =&  L^{-1} (Proj + F \frac{C^{-1}}{2|v_p^*|^2} B),\\
 L& =& \left(
a_j^i \right)_{0\leq i \leq p-1; 1\leq j \leq p}, \; \mbox{{\rm
where $i$ denotes the rows and $j$ the
columns}},\\
F&=&\left(
\begin{array}{cc}
0_p'
 &  \left(\pi_j \; i\; a_j^{i-1}\right)_{i=2,\cdots,p, \;  j=1,\cdots, p}
\end{array}\right)', \\
&& \mbox{{\rm where $i$ denotes the rows, $j$ denotes the
columns,}}\\
Proj && \mbox{{\rm is the projection of $\BBc^{(p+1)^2}$ on $\BBc^p$
such that }} \\
&&\forall w \in \BBc^{(p+1)^2}, \; Proj(w)=v, \;
v=(w_1,\ldots,w_{p})'.\\
\end{eqnarray*}
\end{Corollary}
\begin{Remark} The proof of Theorem \ref{propspeed} is obtained using similar
arguments as in Gassiat {\rm \&} Gautherat's proof of Theorem 4.2
[1999]. The gain, we obtain, with our estimation procedures (without
penalty term) is the asymptotic marginal distribution of both
$(\widehat{\xi}_0 - \xi_0)$ and $(\widehat{\sigma}_{0}-\sigma_0)$.
This is an essential point to obtain the asymptotic distribution of
both $\widehat a$ and $\widehat \Pi$.
\end{Remark}

\section{Simulation study}\label{Simu:Sec}
In this section, the estimates of $\sigma_0$, $\xi_0$, $(a_i)_{1\leq
i \leq p}$ and $(\pi_i)_{1\leq  i \leq p}$ are provided using
 our theoretical estimation procedures. To give stable results,
 each estimation value is an average over $N=100$ independent simulations of
 sequence of $n$ observations. These estimation values are  denoted
 by
 ${\widehat E}$ in the arrays below. The stability of each estimation value is measured by  "std",
 the empirical standard
 deviation calculated over the simulation runs.
The simulations which lead to a negative value of ${\widehat
 \Pi}$ are eliminated and  $N_{elim}$ in the arrays below is the number of those eliminated simulations.

Two models  are considered in this simulation study: the mixture
model and the autoregressive model that is, each simulation sequence
$((Y_{t,l})_{1\leq t \leq n,\; 1 \leq l \leq N})$ is issued of one
of these models which are particular cases of model defined by
(\ref{model}). In both cases, the filter has a finite length, so we
identify the inverse filter $\theta=s(\xi_0)$ to $\xi_0$. We
restrict ourselves to $p=3$ and we deal with two values of
$\sigma_0$: $\sigma_0=0.05$ and $\sigma_0=1$; a small $\sigma_0$
(respectively a large $\sigma_0$) give  small perturbations (resp.
large perturbations) to the corresponding non-noisy model defined by
(\ref{model}) so that it is more complicated to estimate it well in
the case of a large $\sigma_0$. To illustrated the asymptotic
efficiency of our method, we deal with several $n$
($n=50,\;100,\;250,\;500,\;750,\;1000,\;1500,\;2000$) and $k(n)$
($k(n)=1,2$) but only few significant results are presented here. We
use the function "fsolve" in MATLAB version 7 which allows to find a
root of a given function. The problem of the use of such a function
is its very sensitivity to the starting point since it searches a
zero near the starting point. To overcome this difficulty, we try
several starting points in order to select as the initial point the
one which seems more stable during the simulations. To avoid
problems relied on both scale and delay of $\theta$, we fix the
scale and delay in considering ${\widehat \theta}=({\widehat
\theta}_{k})_{|k| \leq k(n)}$ such that $\|{\widehat
\theta}\|_{2,n}=1$ and $({\widehat \theta}_{-k(n)})$ has the
greatest value among all the components of ${\widehat \theta}$.

\begin{figure}[hbtp]
\begin{center}
\includegraphics[width=8cm]{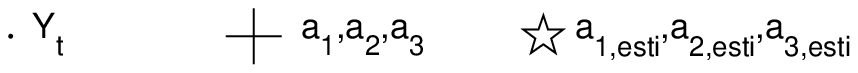}
\\
\includegraphics[width=9cm]{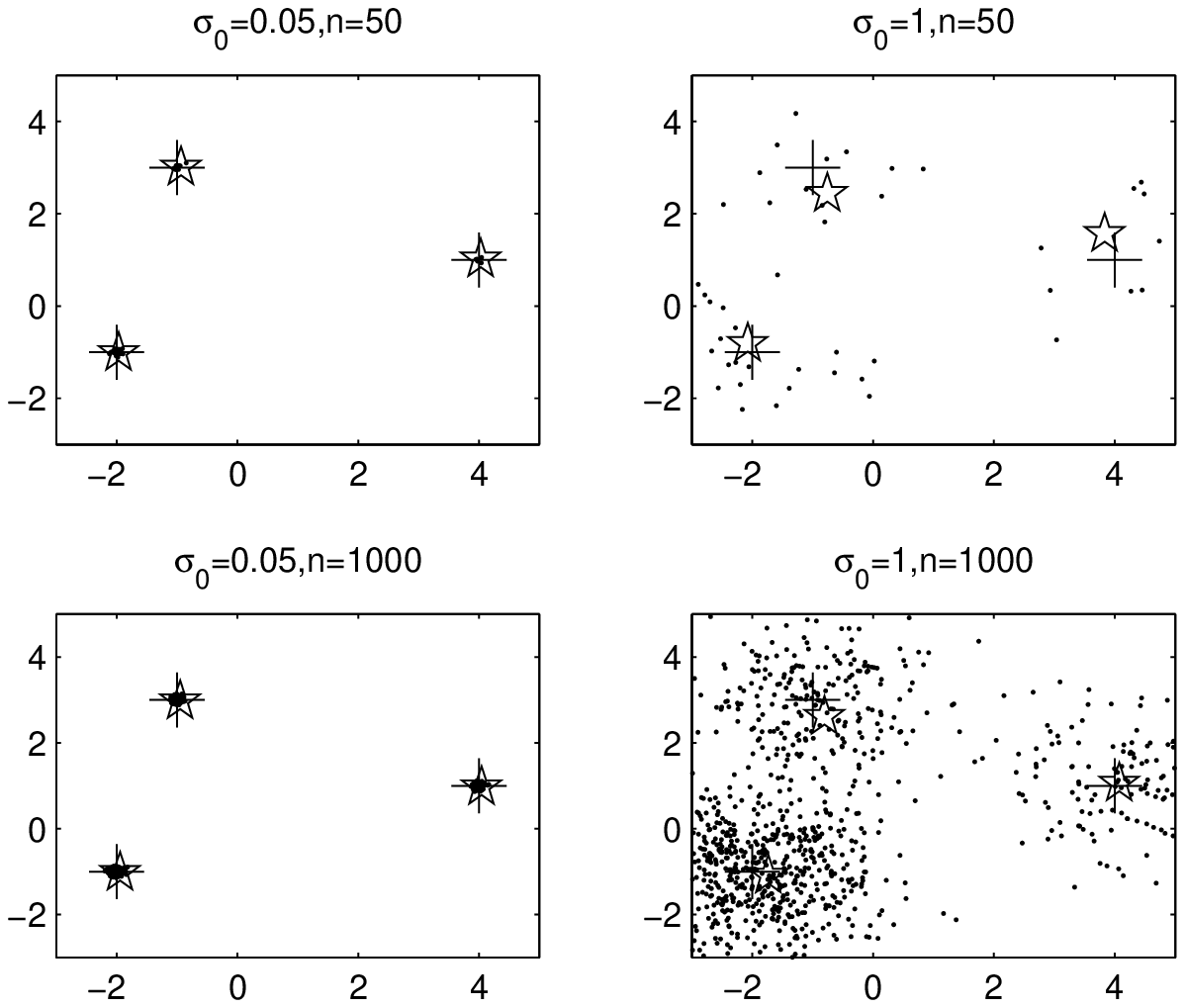}
\end{center}
\caption{{\small Mixture Model, $k(n)=1$. }}
\label{simu-mel-kn1}\end{figure}

\begin{figure}[hbtp]
\begin{center}
\includegraphics[width=8cm]{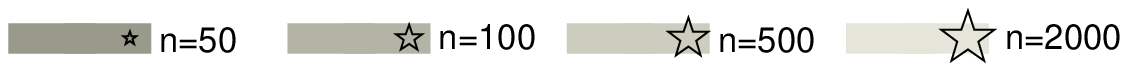}\\
\includegraphics[width=5cm]{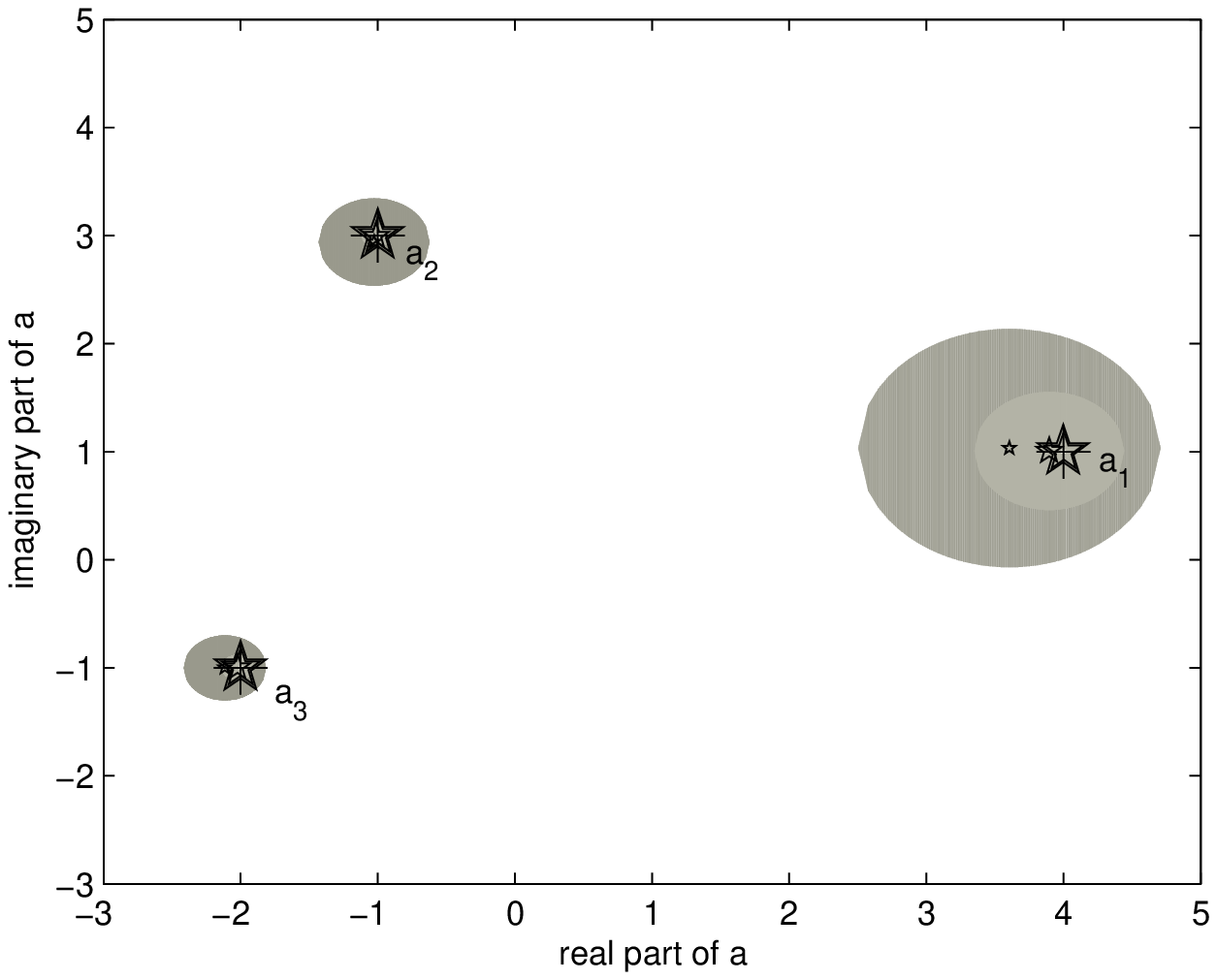}\includegraphics[width=5cm]{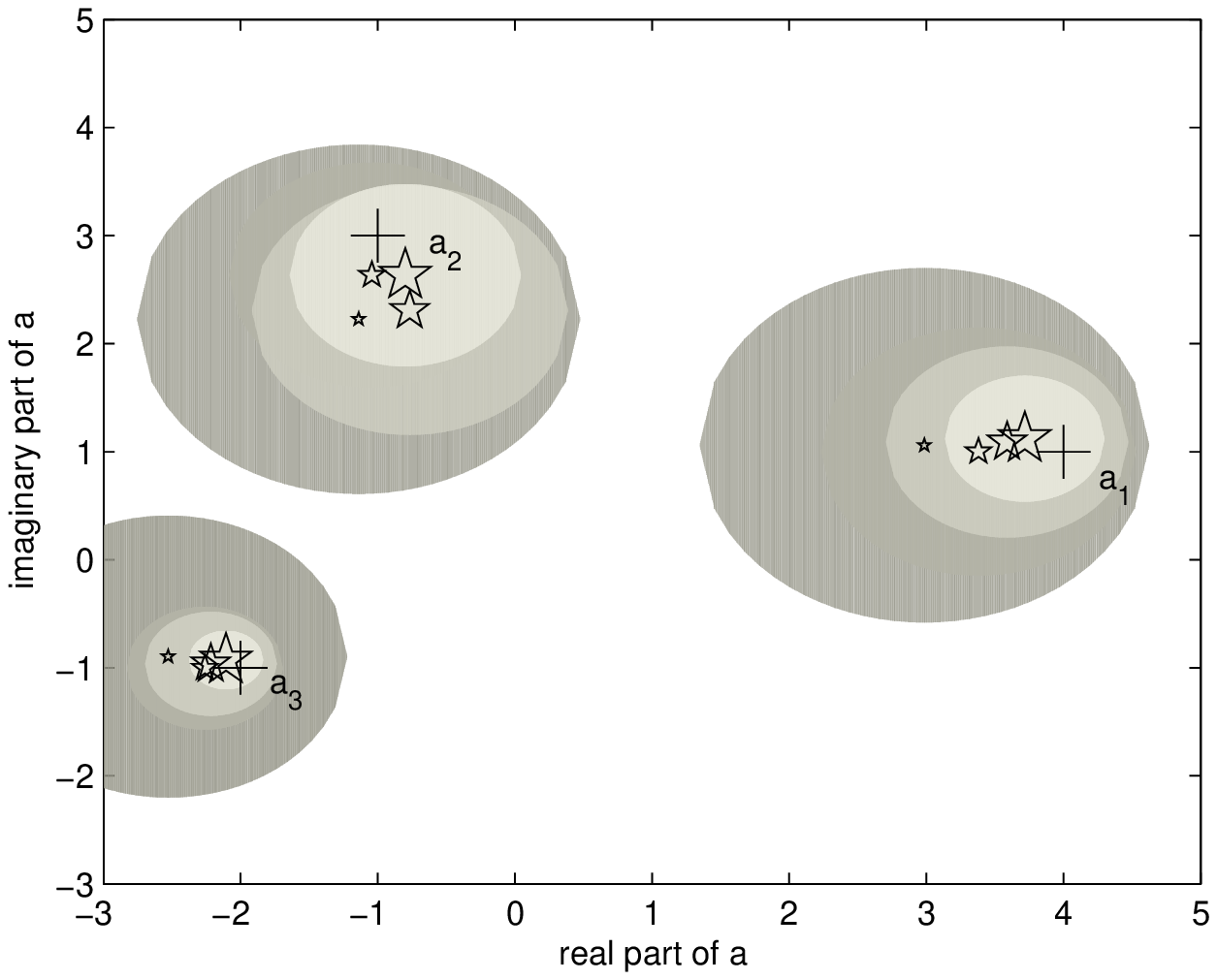}\\
\includegraphics[width=5cm]{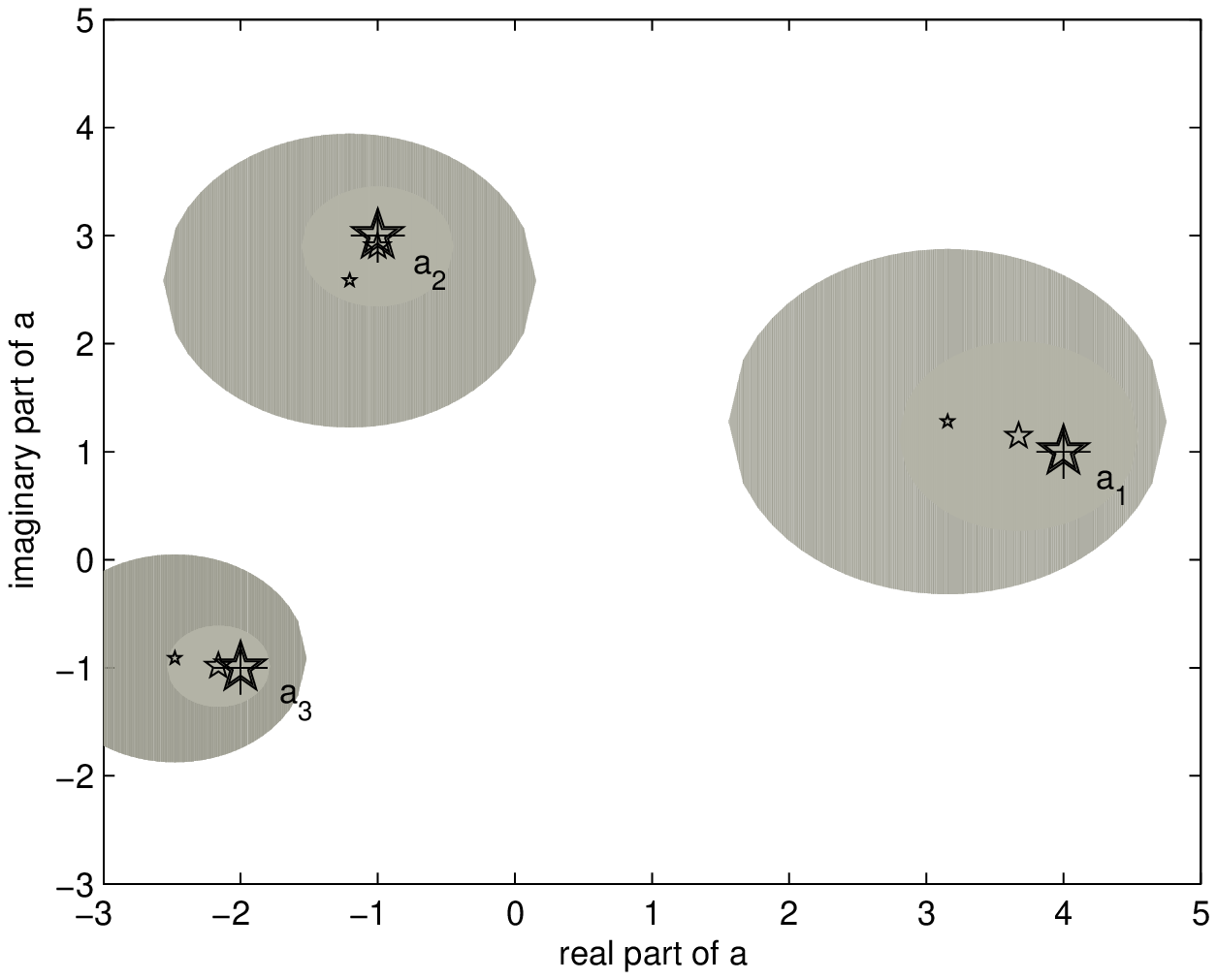}\includegraphics[width=5cm]{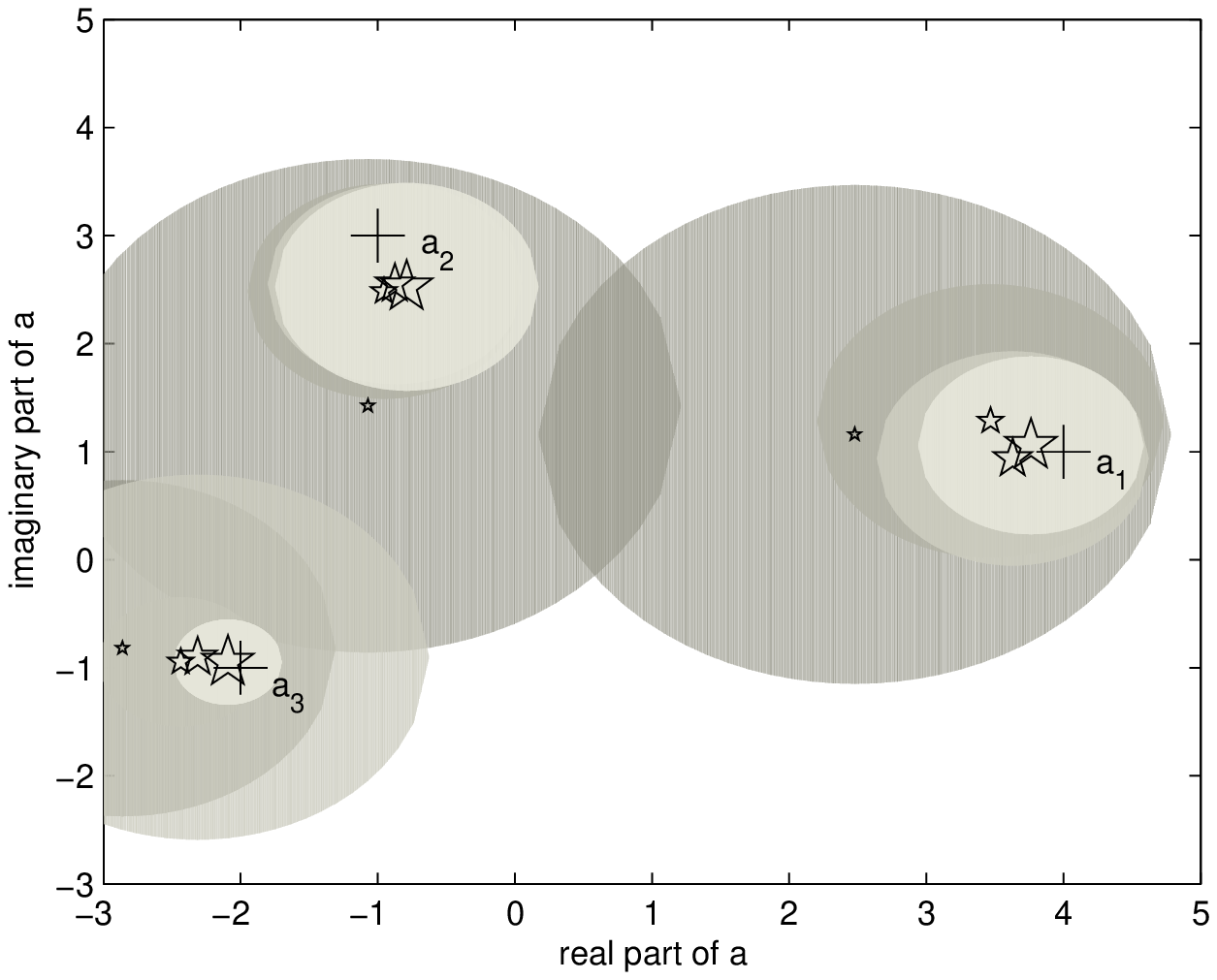}
\caption{ Evolution of ${\widehat a}_1, {\widehat a}_2$ and
${\widehat a}_3$ with $n$. On the first line $k(n)=1$, on the second
line $k(n)=2$. For the first row $\sigma_0=0.05$ and for the second
row $\sigma_0=1$.} \label{figunif-mel-kn1-a}
\end{center}
\end{figure}

\begin{figure}[hbtp]
\begin{center}
\includegraphics[width=8cm]{titresimu.eps}
\\
\includegraphics[width=8cm]{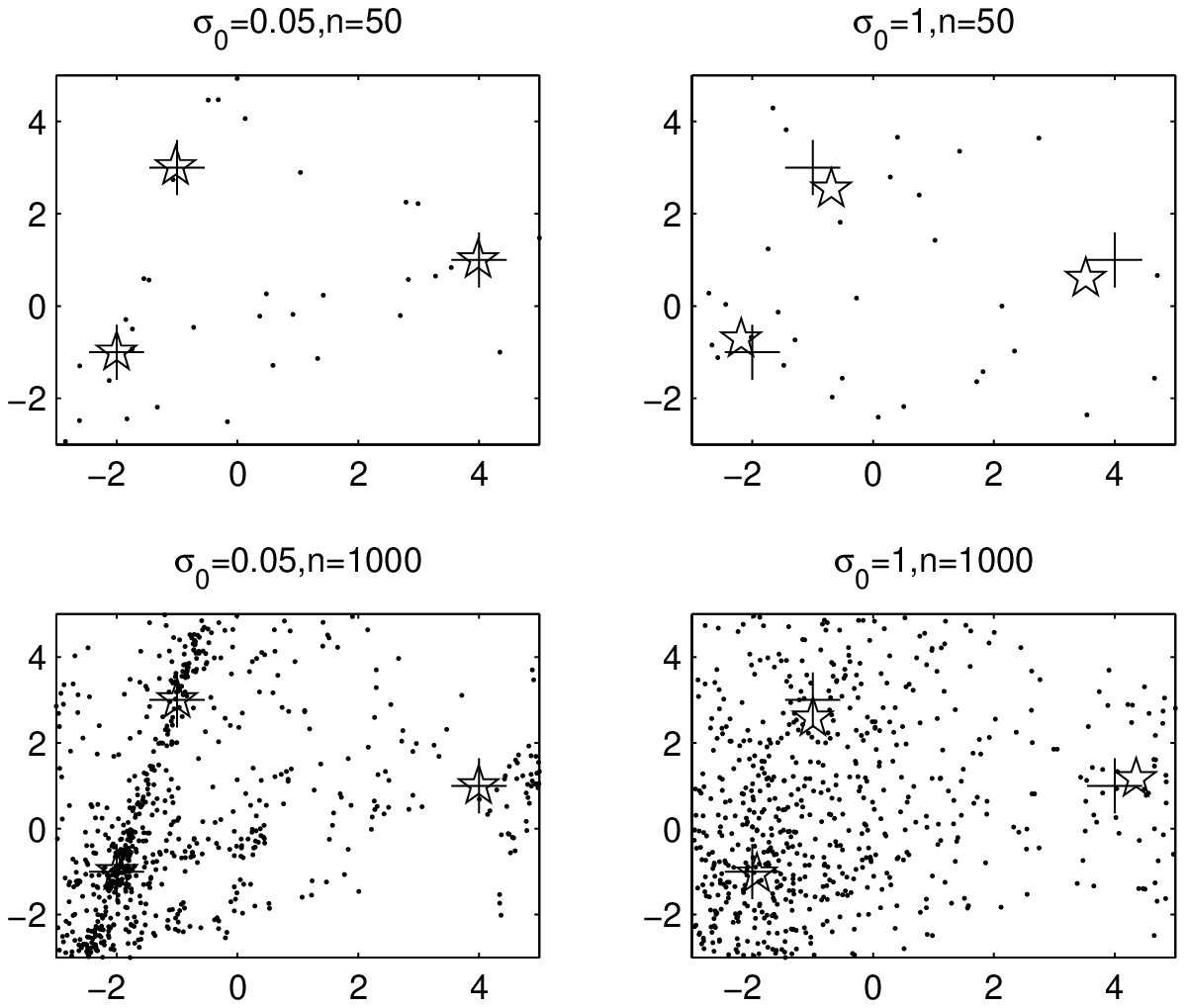}
\caption{{\small Second order autoregressive, $k(n)=1$. }}
\label{simu-ar2-kn1}
\end{center}
\end{figure}

\begin{figure}[hbtp]
\begin{center}
\includegraphics[width=8cm]{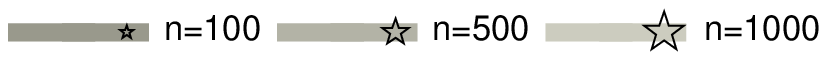}
\\
\includegraphics[width=5cm]{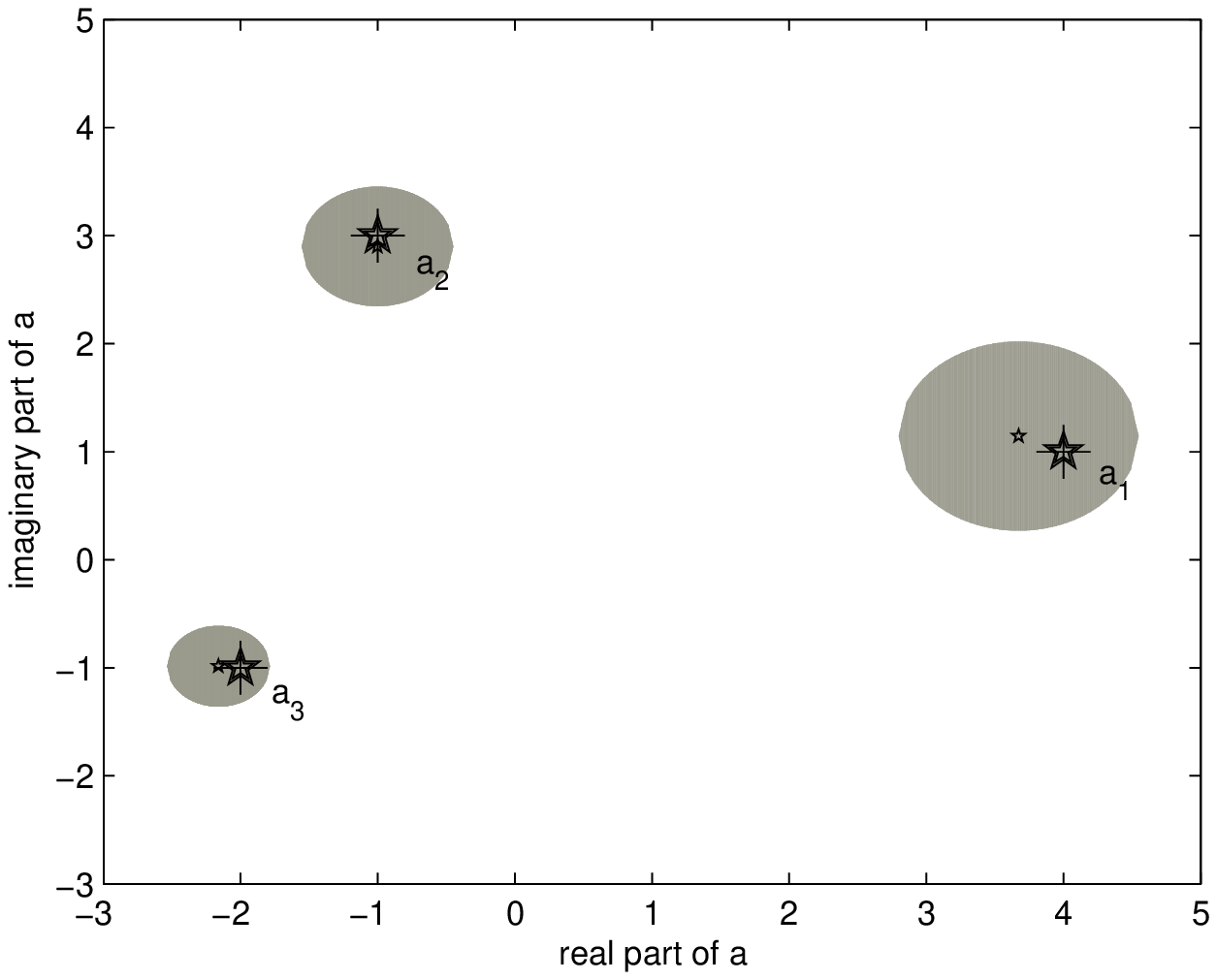}\includegraphics[width=5cm]{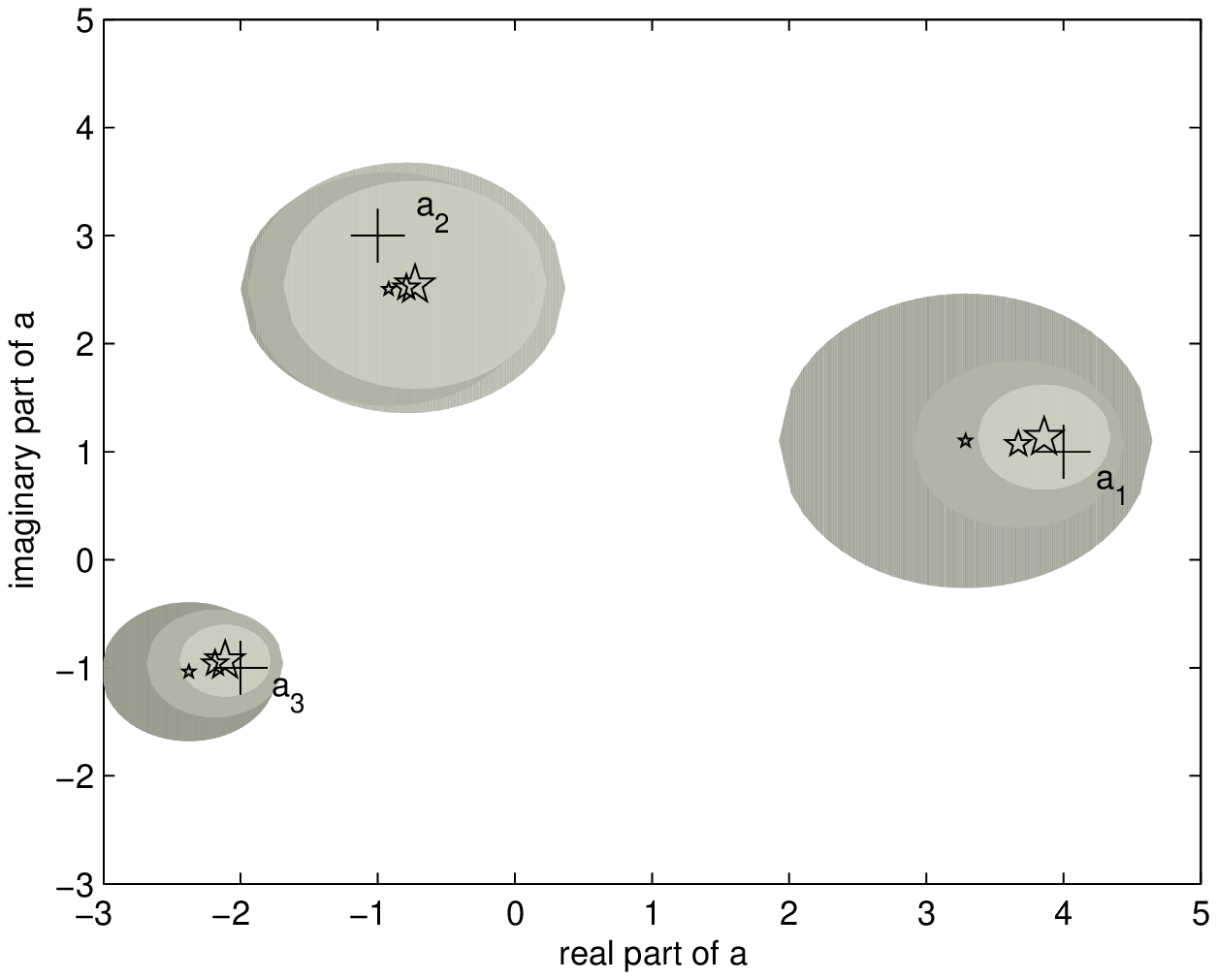}\\
\includegraphics[width=5cm]{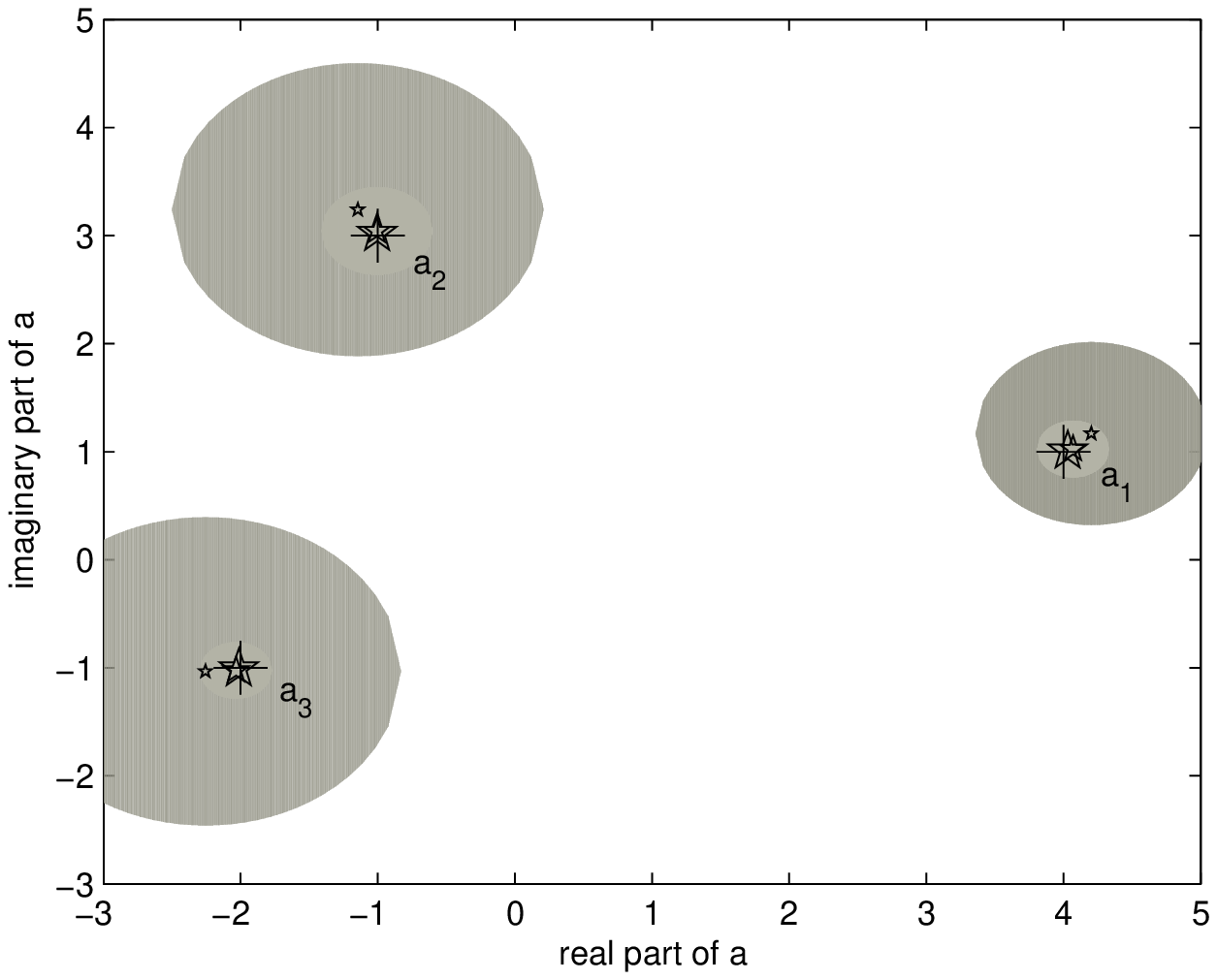}\includegraphics[width=5cm]{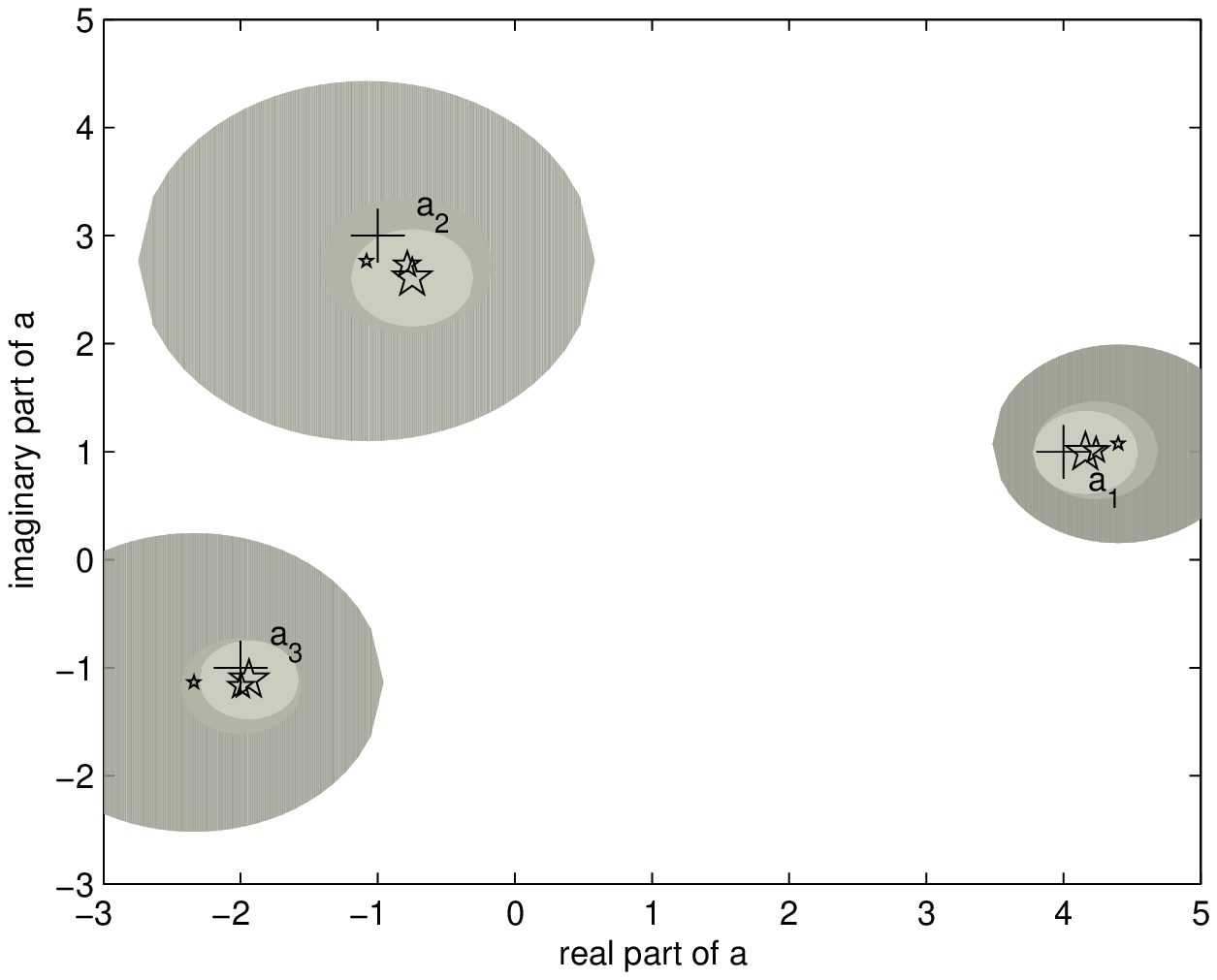}
\caption{ Evolution of ${\widehat a}_1, {\widehat a}_2$ and
${\widehat a}_3$ with $n$. On the first line $k(n)=1$, on the second
line $k(n)=2$. For the first row $\sigma_0=0.05$ and
 for the second row $\sigma_0=1$. } \label{figunif-ar2-a}
 \end{center}
\end{figure}

\noindent {\bf Mixture model.} The observed sequence $Y_1,
\ldots,Y_n$ is given by
$$ Y_t= X_t + \sigma_0 W_t, \; \forall t \in \{1, \ldots, n\},
 $$
so that the filter and its inverse coincide i.e.
$\theta_t=u_t=\delta_0(t) \; \forall t \in \BBz$.
One must note that the model is over-parameterized since we choose
$k(n)=1$ (first array) and $k(n)=2$  (second array) whereas the true
inverse filter is reduced to one value $\theta_0=1$ which
corresponds to  $k(n)=0$. The cases $n=50$ and $n=2000$ are
presented in the arrays below. When $\sigma_0=0.05$, the estimation
values of $(a_i)_{1 \leq i \leq 3}$ and $(\pi_i)_{1 \leq i \leq 3}$
are really good even for a small $n$ ($n=50$) and $k(n)=2$. For both
$\sigma_0=0.05$ and $\sigma_0=1$ and both $k(n)=1$ and $k(n)=2$, the
estimation values of $\sigma_0$ and $\theta$ are strongly different
between $n=50$ and $n=2000$: they are much more better for $n$
large. These arrays illustrate that it is more difficult to obtain a
good estimation when $\sigma$ is large or/and when $k(n)=2$ or/and
when $n=50$.  It is worthwhile to note that for $\sigma_0=1$ the
number of non-used sequences is $N_{elim}=30$ ($n=50$, $k(n)=1$),
$N_{elim}=53$ ($n=2000$, $k(n)=1$), $N_{elim}=38$ ($n=50$, $k(n)=2$)
and $N_{elim}=45$ ($n=2000$, $k(n)=2$), so that it would probably
mean that we do not find an interesting starting point; moreover it
entails a strong variability on the values (see the $std$). Figure
\ref{simu-mel-kn1} represents the $n$ observations $Y_1,\ldots,Y_n$,
the support points $a_1,a_2,a_3$ and their estimates ${\widehat
a}_1$, ${\widehat a}_2$ and ${\widehat a}_3$ which are denoted
$a_{i,\mbox{{\rm esti}}}, \; i=1,2,3$. When $\sigma_0$ is small
($\sigma_0=0.05$), the observations are concentrated on the true
support points and the estimation for $k(n)=1$ is visually very good even for
a small sample (Figure \ref{simu-mel-kn1}, left side). When $\sigma_0$ is large
($\sigma_0=1$), the observations are more scattered over the square
even one can distinguish three attractive areas where the support
points lie. For this case and $k(n)=1$, the estimation with $1000$ observations
is visually much more better than the one with $50$ observations (Figure \ref{simu-mel-kn1}, right side).
The same phenomenon is observed in Figure \ref{figunif-mel-kn1-a}
but in addition Figure \ref{figunif-mel-kn1-a} illustrates the
improvement of the estimation: when $n$ is increasing, the estimates
approach the true values and their empirical standard deviations
tend to zero (except for one case).
\vspace{0.3cm} \noindent{\bf Second order autoregressive model.} The
observed sequence $Y_1,\ldots,Y_n$ is given by
$$ Y_t= \tilde{Y}_t + \sigma_0 W_t,  \; \forall t \in \{1, \ldots,
n\},$$
where
$$\tilde{Y}_t =  \sum_{k=-\infty }^0 u_k X_{t-k}
\Longleftrightarrow
 X_t= \sum_{k=0}^2 \theta_k \; \tilde{Y}_{t-k},$$
with $(\theta_0,\theta_1,\theta_2)=(\frac 67,-\frac 27,\frac
37)=(0.8571, -0.2857, 0.4286)$.
One must note that when we choose $k(n)=2$, the model is
over-parameterized since the true model corresponds to $k(n)=1$. In
the left hand side of the third table ($k(n)=1$) with a small
$\sigma_0=0.05$, the estimations of all the parameters are quite
good and they are quite similar for both large $n$ and small $n$.
All other cases i.e $\sigma_0=1$ with $k(n)=1$ and $k(n)=2$, large
$n$ ($n=1000$ or $n=2000$) lead to an improvement in the estimation
values since some estimations for $n=50$ are very far from the true
values (see both the right hand side of the third table and the
fourth table). However, in some cases it seems that we do not take
an interesting starting point since the variability of the results
is too large (see for example the $std$ for $k(n)=1$ and
$\sigma_0=1$) combined with an important number of eliminated
simulations (see $N_{elim}$ for $k(n)=1$ and $\sigma_0=1$). One can
see in Figure \ref{simu-ar2-kn1} that the observations are more
dispersed over the square than those of the mixture model (Figure
\ref{simu-mel-kn1}). It would mean that an autoregressive model is
more difficult to estimate than a mixture model. When $\sigma_0$ is
small ($\sigma_0=0.05$), Figure \ref{simu-ar2-kn1} and Figure
\ref{figunif-ar2-a} show that the support's points $(a_i)_{1 \leq i
\leq 3}$ are well estimated
 for $k(n)=1$ even for a small $n$; this is not the case for
$k(n)=2$ since the estimates seem visually (see Figure
\ref{figunif-ar2-a}) far from the true values even for large $n$
($n=1000$).
\medskip

\begin{figure}[hbtp]
\begin{center}
\includegraphics[width=5cm]{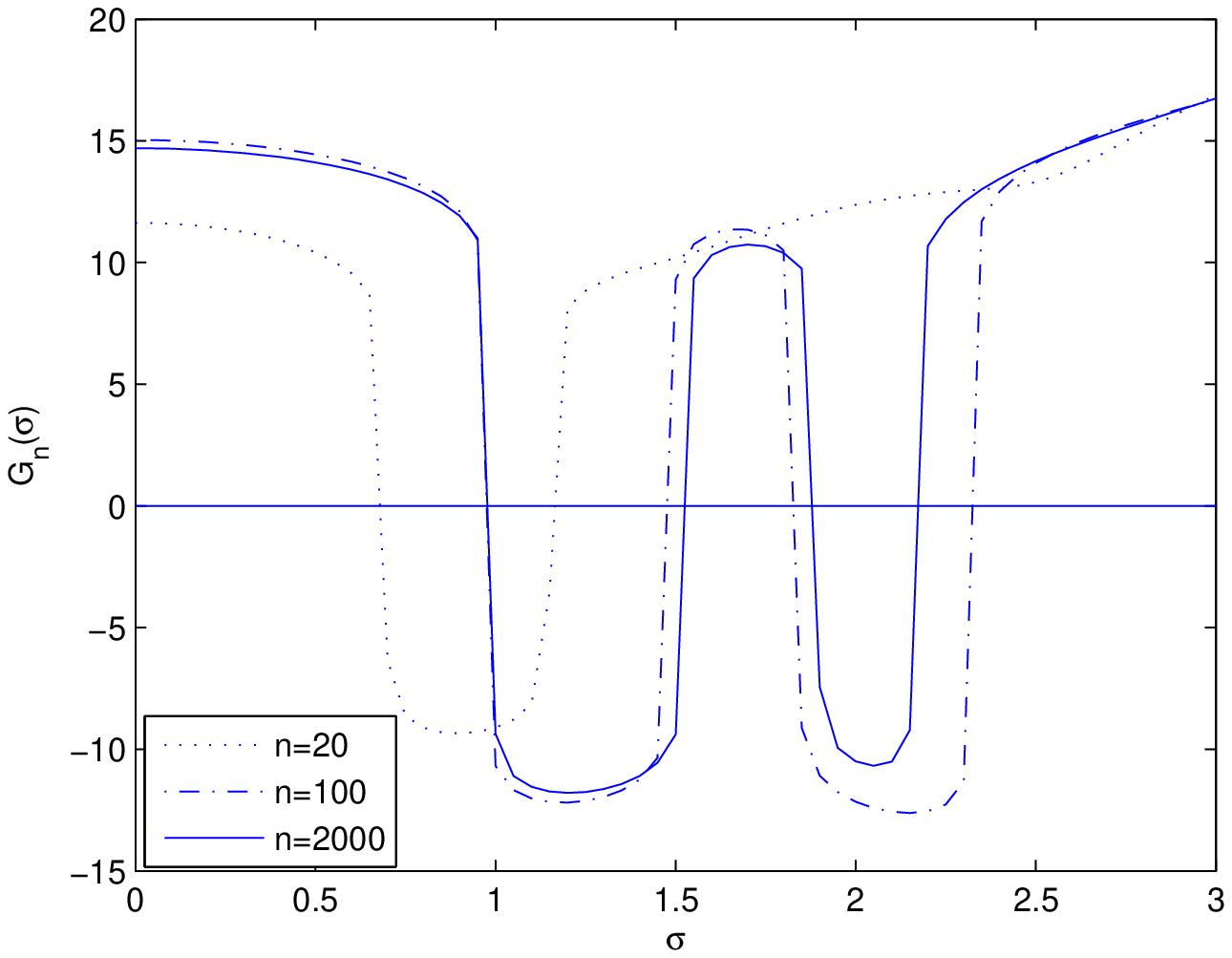}
\includegraphics[width=5cm]{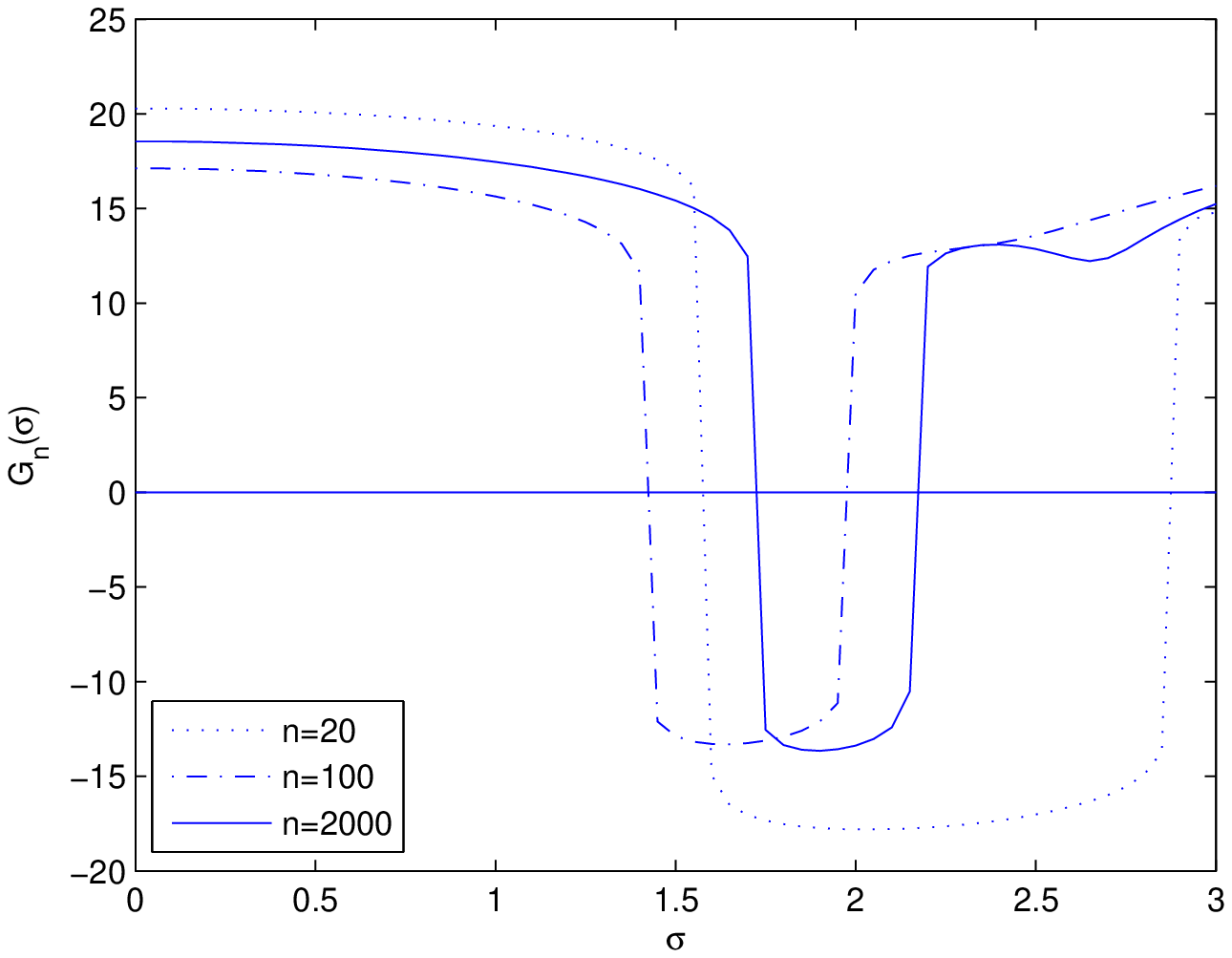}
\end{center}
\medskip
\caption{Graph of $G_{n,\bar s (\xi)}$ for $k(n)=2$. Second order
autoregressive, $\sigma_0=1$.  } \label{figcoupe}
\end{figure}

\noindent{\bf  Numerical illustration of the existence of the
estimate $\widehat{\sigma}_{0}$.} Let us consider the second order
autoregressive model with $\sigma_0=1$ and  let  $k(n)$ be equal to
$2$. This corresponds to the most difficult estimation problem as it
is illustrated in the above simulations. Figure \ref{figcoupe}
represents the graph of  the function $\sigma \rightarrow G_{n,{\bar
s}(\xi)}(\sigma)=sign(J_n(\sigma,\xi))*log(|J_n(\sigma,\xi)|+1)$
which has the same roots as the function $\sigma \rightarrow
J_n(\sigma,\xi)$. In the left hand side, we consider the particular
case of  ${\bar s}(\xi_0)=(\theta_0,\theta_1,\theta_2,0,0)$ which is
the true value of the inverse filter, whereas in the right hand
side, we consider the particular case of $\bar s(\xi)=(1,0,0,0,0)$
which differs from the true value of the inverse filter (this filter
corresponds to a mixture model). In both cases, one can see that
$G_{n,{\bar s}(\xi)}$ admits some zeros. In the left hand side, one
must note that the convergence is achieved very quickly and
accurately.

\medskip

\noindent{\bf Importance of the choice of the starting point in the
algorithm.} One can see in Figure \ref{figcoupe} that the first zero
of the function $J_n$ is  neither  over $\sigma_0$ nor under
$\sigma_0$. Actually, the Matlab toolbox algorithm searches the zero
of the function which is the nearest to the starting point. This
starting point being a $(2k(n) +2)$-dimensional vector, the
algorithm searches a zero in $(2 k(n) +2)$ directions. So an
iterative stochastic algorithm, which is able to find a
multidimensional zero, would be an useful tool since the gradient of
$J_n$ could be formally computed.

\begin{table}
\begin{center}
\begin{tiny}
\begin{tabular}{|c|c|rl|}
\hline $\sigma_0=0.05 $& $n$ & $\widehat{E}$ &$\pm std$
\\
$k(n)=1$ & & &  \\
 \hline $\widehat{\sigma}_{0}$ & 50&0.1397& $\pm 0.2125$
\\
& 2000&0.0554  &     $\pm               0.0114$
\\
\hline $\widehat{\theta}_0$ & 50&0.9627   &$\pm 0.0962$
\\
& 2000& 0.9999                &  $ \pm                 0.0002$
\\
\hline $\widehat{\theta}_1 $& 50&0.0732  &$ \pm 0.1866$
\\
& 2000&0.0000  &  $\pm                  0.0128$
\\
\hline $\widehat{\theta}_2$ & 50&0.0510 & $\pm 0.1479$
\\
& 2000&0.0006 & $\pm                     0.0063$
\\
\hline $\widehat{a}_3 $& 50&-2.1146 - 1.0012 i& $\pm 0.2996$
\\
& 2000&-2.0004  - 0.9996 i& $ \pm 0.0130$
\\
\hline $\widehat{a}_2 $& 50&-1.0265 +2.9403 i &$ \pm  0.4038$
\\
& 2000&-1.0002 + 3.0000 i&$ \pm 0.0129$
\\
\hline$ \widehat{a}_1$ & 50&3.6058 + 1.0336 i & $\pm 1.1026$
\\
& 2000& 3.9984 + 1.0000 i & $ \pm 0.0152$
\\
\hline$ \widehat{\Pi}_1$ & 50&0.5906  &$ \pm 0.0690$
\\
& 2000&0.6000          & $ \pm                 0.0104$
\\
\hline $\widehat{\Pi}_2$ & 50&0.2468  & $\pm 0.0671$
\\
& 2000&0.2509    &  $\pm                   0.0091$
\\
\hline $\widehat{\Pi}_3$ & 50&0.1625  & $\pm 0.0597$
\\
& 2000&0.1492          &   $\pm          0.0078$
\\
\hline $N_{elim}$ & 50&4&
\\
& 2000&0&
\\
\hline
\end{tabular}
\begin{tabular}{|c|c|rl|}
\hline
   $\sigma_0=1$& $n$&$ \widehat{E}$& $\pm std $
\\
$k(n)=1$ & & & \\
\hline $\widehat{\sigma}_{0}$ & 50&1.0139  &$ \pm 0.1672$
\\
& 2000&  1.1403   &    $ \pm  0.1936$
\\
\hline $\widehat{\theta}_0$ & 50&0.8821   &$\pm 0.1100$
\\
& 2000&   0.9107                &  $ \pm 0.1296$
\\
\hline $\widehat{\theta}_1 $& 50&0.1192      & $\pm 0.3169$
\\
& 2000& 0.0931       &  $\pm  0.2834$
\\
\hline $\widehat{\theta}_2 $& 50&0.1087  & $\pm  0.2935$
\\
& 2000& 0.0489      &$ \pm   0.2566$
\\
\hline $\widehat{a}_3 $& 50& -2.5278-  0.8964 i& $\pm 1.3044$
\\
& 2000& -2.1064   -   0.9275 i&  $\pm 0.2689$
\\
\hline$ \widehat{a}_2$& 50&-1.1398 + 2.2255  i& $\pm 1.6165$
\\
& 2000&-0.7987  + 2.6327  i& $\pm 0.8427$
\\
\hline$ \widehat{a}_1 $& 50&2.9865+  1.0596 i &$ \pm 1.6399$
\\
& 2000&3.7181  +1.1201 i & $ \pm  0.5812$
\\
\hline $\widehat{\Pi}_1$ & 50& 0.4853  & $\pm 0.1214$
\\
& 2000& 0.5681         &$  \pm 0.0865$
\\
\hline$ \widehat{\Pi}_2 $& 50&  0.3374             & $\pm 0.1397$
\\
& 2000&  0.2923            & $ \pm    0.0939$
\\
\hline$ \widehat{\Pi}_3$ & 50& 0.1773            & $\pm  0.0686$
\\
& 2000&  0.1396            &  $ \pm  0.0233$
\\
\hline $N_{elim}$ & 50&30&
\\
& 2000&53&
\\
\hline
\end{tabular}
\end{tiny}
\caption{
$(\theta_0,\theta_1,\theta_2,a_1,a_2,a_3,\pi_1,\pi_2,\pi_3)=
 (1,0,0,4+i,-1+3i,-2-i,0.6,0.25,0.15)$. The Starting point of left
array is
$(\sigma_{ini},\theta_{ini,0},\theta_{ini,1},\theta_{ini,2})=(
0.001,1.2,-0.5,0.02)$  and the starting point of right array is
$(\sigma_{ini},\theta_{ini,0},\theta_{ini,1},\theta_{ini,2})=(
0.1,1.2,-0.4,0.2)$.}\label{tablemixturea}
\end{center} \end{table}

\begin{table}
\begin{center}
\begin{tiny}
\begin{tabular}{|c|c|rl|}
\hline
  $\sigma_0=0.05$& $n$&$ \widehat{E}$&$\pm std$
\\
$k(n)=2$ & & & \\
\hline $\widehat{\sigma}_{0} $& 50&0.3913& $\pm 0.3550$
\\
& 2000& 0.0542    &     $\pm               0.0078$
\\
\hline $\widehat{\theta}_0$ & 50&0.8362    &$\pm 0.1686$
\\
& 2000&0.9999       &  $ \pm      0.0001$
\\
\hline$ \widehat{\theta}_1$ & 50& 0.0674 & $\pm 0.2417$
\\
& 2000&-0.0003   & $ \pm          0.0053$
\\
\hline $\widehat{\theta}_2$ & 50&0.0840   &$ \pm  0.2086$
\\
& 2000&0.0001         & $\pm    0.0080$
\\
\hline $\widehat{\theta}_3$ & 50& 0.1308 & $\pm 0.2520$
\\
& 2000&0.0000            & $\pm  0.0047$
\\
\hline $\widehat{\theta}_4$ & 50& 0.1026     & $\pm 0.2664$
\\
& 2000&-0.0000      &$ \pm       0.0054$
\\
\hline $ \widehat{a}_3 $& 50&-2.4801 -0.9127  i &$\pm 0.9610$
\\
& 2000&-1.9996 - 1.0002 i  & $ \pm  0.0112$
\\
\hline $\widehat{a}_2 $& 50&-1.2051 + 2.5834  i    &$ \pm 1.3587$
\\
& 2000&-0.9997 + 3.0000 i & $\pm  0.0115$
\\
\hline $\widehat{a}_1$ & 50&3.1553 +  1.2781 i&$ \pm 1.5967$
\\
& 2000&3.9999 +  1.0001 i &  $\pm0.0109$
\\
\hline $\widehat{\Pi}_1$ & 50& 0.5211  & $\pm 0.1270$
\\
& 2000&0.5997  &  $\pm           0.0115$
\\
\hline $\widehat{\Pi}_2$ & 50& 0.3011   & $\pm 0.1314$
\\
& 2000&0.2510     &  $\pm        0.0091$
\\
\hline $\widehat{\Pi}_3$ & 50&0.1778  & $\pm 0.0812$
\\
& 2000& 0.1493          &   $\pm      0.0084$
\\
\hline$ N_{elim}$ & 50&20&
\\
& 2000&0&
\\
\hline
\end{tabular}
\begin{tabular}{|c|c|rl|}
\hline $\sigma_0=1$& $n$& $\widehat{E}$&$\pm std$
\\
$k(n)=2$ & & & \\
\hline $\widehat{\sigma}_{0} $& 50& 1.0467  & $\pm 0.1906$
\\
& 2000&1.1649      &     $\pm               0.2807$
\\
\hline $\widehat{\theta}_0$ & 50& 0.7038     &$\pm 0.1283$
\\
& 2000&0.9195                    &   $\pm               0.1405$
\\
\hline $\widehat{\theta}_1$ & 50&0.1393 & $\pm 0.3311$
\\
& 2000& -0.0087             & $ \pm                    0.1796$
\\
\hline $\widehat{\theta}_2$ & 50& 0.1962  & $\pm  0.2900$
\\
& 2000& 0.0232      &$ \pm                    0.1728$
\\
\hline $\widehat{\theta}_3$ & 50&0.1724  &$ \pm 0.2707$
\\
& 2000&  0.0077   &$ \pm                      0.2007$
\\
\hline$ \widehat{\theta}_4$ & 50&0.0998           & $\pm 0.3609$
\\
& 2000&   0.0077         & $\pm                      0.1858$
\\
\hline $ \widehat{a}_3$ & 50&-2.8623- 0.8196 i & $\pm 1.5551$
\\
& 2000& -2.0936 -0.9478 i &  $\pm   0.3934$
\\
\hline $\widehat{a}_2$ & 50&-1.0714 + 1.4250  i    & $\pm 2.2831$
\\
& 2000&-0.7886 + 2.5256 i& $\pm 0.9606$
\\
\hline $\widehat{a}_1 $& 50&  2.4786 +  1.1597 i & $\pm 2.3074$
\\
& 2000& 3.7623 +  1.0584 i & $ \pm  0.8220$
\\
\hline$ \widehat{\Pi}_1$ & 50&0.4341   & $\pm 0.1474$
\\
& 2000&0.5436     &  $\pm                 0.1302$
\\
\hline$ \widehat{\Pi}_2 $& 50&0.3476     & $\pm 0.1407$
\\
&2000& 0.3161                &  $\pm                  0.1294$
\\
\hline $\widehat{\Pi}_3 $& 50&0.2183            & $\pm  0.1353$
\\
& 2000&  0.1403             &   $\pm              0.0317$
\\
\hline $N_{elim}$ & 50&38&
\\
& 2000&45&
\\
\hline
\end{tabular}
\end{tiny}
 \caption{  $(1,0,0,4+i,-1+3i,-2-i,
0.6,0.25,0.15)$. Starting point of left array
$(\sigma_{ini},(\theta_{ini,i})_{0\leq i \leq 4})\!\!\!=\!\!\!(
0.01,0.3,-0.1,0.2,-0.1,0.1)$ and Starting point right array
$(\sigma_{ini},(\theta_{ini,i})_{0\leq i \leq 4})\!\!\!=\!\!\!(
0.1,0.3,-0.1,0.2,-0.1,0.1)$.}\label{tablemixtureb}
\end{center} \end{table}
\newpage
\begin{table}
\begin{center}
\begin{tiny}
\begin{tabular}{|c|c|rl|}
\hline $\sigma_0=0.05$& $n$&$\widehat{E}$&$\pm std$
\\
$k(n)=1 $& & &  \\
\hline $\widehat{\sigma}_{0}$ & 50&0.0565& $\pm 0.0133$
\\
& 1000&  0.0591              &     $\pm               0.0161$
\\
\hline $\widehat{\theta}_0 $& 50&0.8627   &$\pm  0.0045$
\\
& 1000&0.8636               &  $ \pm                   0.0068$
\\
\hline $\widehat{\theta}_1$ & 50&-0.2786 &$ \pm 0.0076$
\\
& 1000& -0.2776           &  $\pm                   0.0098$
\\
\hline $\widehat{\theta}_2$ & 50&0.4220  & $\pm 0.0051$
\\
& 1000& 0.4207     &$ \pm                      0.0078$
\\
\hline $\widehat{a}_3$ & 50&-2.0130 - 1.0067i&$\pm  0.0182$
\\
& 1000&-2.0174 - 1.0095i  & $ \pm  0.0240$
\\
\hline $\widehat{a}_2$ & 50&-1.0090 + 3.0197i&$ \pm  0.0276$
\\
& 1000& -1.0088 + 3.0246i  & $\pm  0.0319$
\\
\hline $\widehat{a}_1$ & 50&4.0261 + 1.0081i & $\pm  0.0299$
\\
& 1000& 4.0338 + 1.0082i  &  $\pm  0.0395$
\\
\hline $\widehat{\Pi}_1 $& 50&0.6087  &$ \pm  0.0539$
\\
& 1000& 0.5998    & $ \pm                  0.0147$
\\
\hline $\widehat{\Pi}_2$ & 50&0.2517         & $\pm  0.0489$
\\
& 1000& 0.2504   &  $\pm                  0.0133$
\\
\hline $\widehat{\Pi}_3 $& 50&0.1396    & $\pm  0.0326$
\\
& 1000&0.1498            &   $\pm          0.0105$
\\
\hline $N_{elim}$ & 50&0&
\\
& 1000&0&
\\
\hline
\end{tabular}
\begin{tabular}{|c|c|rl|}
\hline $   \sigma_0=1$&$ n$&$ \widehat{E}$&$\pm std$
\\
$k(n)=1 $& & &   \\
\hline $\widehat{\sigma}_{0}$ & 100&1.1764& $\pm 0.2959$
\\
& 1000& 1.2805 &     $\pm 0.3006$
\\
\hline $\widehat{\theta}_0$ & 100&0.8263  & $\pm 0.0572$
\\
& 1000&   0.7814             &   $\pm 0.0541$
\\
\hline $\widehat{\theta}_1 $& 100&0.0325 &$ \pm 0.3910$
\\
& 1000& -0.1880          & $ \pm  0.2438$
\\
\hline $\widehat{\theta}_2 $& 100&0.1640  & $\pm 0.3700$
\\
& 1000&  0.4816   & $\pm 0.2503$
\\
\hline$ \widehat{a}_3$ & 100&-3.4094 - 0.8583i   &$\pm  3.4494$
\\
& 1000& -2.3470 - 1.2209i  & $ \pm 1.1879$
\\
\hline $\widehat{a}_2$ & 100&-1.9715 + 3.5025i  & $\pm  3.6024$
\\
& 1000&-1.0180 + 2.7762i &$ \pm 2.0138$
\\
\hline $\widehat{a}_1$ & 100&5.3438 + 1.1887i   & $\pm 5.2553$
\\
& 1000& 3.4259 +  1.0842 i& $ \pm  1.8686$
\\
\hline $\widehat{\Pi}_1 $& 100&0.5078     & $\pm  0.1227$
\\
& 1000& 0.5235      &  $\pm 0.1352$
\\
\hline $\widehat{\Pi}_2 $& 100&0.3075          & $\pm 0.0898$
\\
& 1000&0.2877       &  $\pm 0.0833$
\\
\hline $\widehat{\Pi}_3 $& 100&0.1847       & $\pm 0.0635$
\\
& 1000& 0.1887       &   $\pm 0.1088$
\\
\hline $N_{elim}$ & 100&58&
\\
& 1000&56&
\\
\hline
\end{tabular} \end{tiny} %
\caption{
$(\theta_0,\theta_1,\theta_2,a_1,a_2,a_3,\Pi_1,\Pi_2,\Pi_3)=(6/7,-2/7,3/7,4+i,-1+3i,-2-i,
0.6,0.25,0.15)$. The starting point of left array is
$(\sigma_{ini},\theta_{ini,0},\theta_{ini,1},\theta_{ini,2})=(
0.001,0.5,-0.5,0.5)$ whereas the right array's is
$(\sigma_{ini},\theta_{ini,0},\theta_{ini,1},\theta_{ini,2})=(
0.5,0.6,-0.2,0.2)$.}\label{tableAR2a} \end{center}
\end{table}

\begin{table}
\begin{center}
\begin{tiny}
\begin{tabular}{|c|c|rl|}
\hline $\sigma_0=0.05$& $n$& $\widehat{E}$&$\pm std$
\\
$k(n)=2$ & & & \\  \hline $\widehat{\sigma}_{0}$ & 50& 0.1642 & $\pm
0.3229$
\\
& 1000& 0.0580         & $\pm     0.0080$
\\
\hline $\widehat{\theta}_0$ & 50&0.8090 &$\pm 0.0770$
\\
& 1000&0.8632     &$ \pm         0.0010$
\\
\hline $\widehat{\theta}_1$ & 50& -0.2855      & $\pm 0.0695$
\\
& 1000&-0.2787         &  $\pm    0.0019$
\\
\hline $\widehat{\theta}_2$ & 50&  0.3815   & $\pm 0.0754$
\\
& 1000&  0.4209   & $\pm          0.0011$
\\
\hline$ \widehat{\theta}_3$ & 50& 0.1552 & $\pm  0.1866$
\\
& 1000& 0.0002         & $\pm       0.0018$
\\
\hline $\widehat{\theta}_4$ & 50& -0.1278  &$ \pm 0.1658$
\\
& 1000&-0.0014   &$ \pm          0.0014$
\\
\hline $ \widehat{a}_3$ & 50& -2.5866 - 0.6725i  &$\pm 2.5778$
\\
& 1000&-2.0125 -  1.0080 i  & $ \pm 0.0050$
\\
\hline$ \widehat{a}_2$ & 50&  -1.4271 +  2.8472 i& $\pm  2.0162$
\\
& 1000&-1.0046 +  3.0204 i&$ \pm 0.0068$
\\
\hline$ \widehat{a}_1$ & 50&  4.2549 +   1.4208 i  &$ \pm1.2687$
\\
& 1000& 4.0316 + 1.0070  i& $ \pm  0.0085$
\\
\hline $\widehat{\Pi}_1$ & 50&  0.5268   & $\pm           0.1337$
\\
& 1000&0.6000        &$  \pm     0.0145$
\\
\hline$ \widehat{\Pi}_2 $& 50&  0.2983    & $\pm           0.1169$
\\
& 1000&0.2508           &  $\pm   0.0117$
\\
\hline $\widehat{\Pi}_3$ & 50&  0.1749   & $\pm 0.0664$
\\
& 1000&0.1492     & $  \pm              0.0103$
\\
\hline $N_{elim}$ & 50&3&
\\
& 1000&0&
\\
\hline
\end{tabular}
\begin{tabular}{|c|c|rl|}
\hline $\sigma_0=1$& $n$& $\widehat{E}$&$\pm std$
\\
$k(n)=2$ & & & \\
\hline $\widehat{\sigma}_{0}$ & 100&1.1082 & $\pm 0.1266$
\\
& 2000&1.2655       &    $ \pm      0.0751$
\\
\hline $\widehat{\theta}_0$ & 100&0.8314    &$\pm   0.0613$
\\
& 2000& 0.8390       & $\pm 0.0475$
\\
\hline $\widehat{\theta}_1$ & 100&-0.3027  & $\pm   0.1049$
\\
& 2000& -0.2961  & $\pm  0.1154$
\\
\hline $\widehat{\theta}_2$ & 100& 0.2419   &$ \pm   0.1115$
\\
& 2000& 0.2415          & $\pm   0.1018$
\\
\hline$ \widehat{\theta}_3 $& 100& 0.1408   & $\pm   0.1808$
\\
& 2000& 0.1062         & $\pm 0.1910$
\\
\hline $\widehat{\theta}_4$ & 100&-0.2608    & $\pm   0.1080$
\\
& 2000& -0.2639 &$ \pm 0.0874$
\\
\hline $ \widehat{a}_3 $& 100&-2.3404 - 1.1354i   &$\pm   1.3817$
\\
& 2000& -1.9004 -  1.1264 i &  $\pm 0.2632$
\\
\hline $\widehat{a}_2$ & 100&-1.0821 + 2.7644i      & $\pm 1.6653$
\\
& 2000& -0.7037 + 2.5984 i  & $\pm 0.3441$
\\
\hline $\widehat{a}_1$ & 100&   4.3995 + 1.0729  i   & $\pm 0.9172$
\\
& 2000&   4.1748 +  0.9867 i  & $ \pm 0.3380$
\\
\hline $\widehat{\Pi}_1$ & 100&0.5133   & $\pm     0.1108$
\\
& 2000&  0.5678          &  $\pm    0.0451$
\\
\hline $\widehat{\Pi}_2$ & 100& 0.3115       & $\pm   0.0991$
\\
& 2000&  0.2799            & $\pm 0.0363$
\\
\hline $\widehat{\Pi}_3 $& 100&0.1752         & $\pm   0.0544$
\\
& 2000&  0.1524      & $\pm 0.0142$
\\
\hline $N_{elim}$ & 100&11&
\\
& 2000&14&
\\
\hline
\end{tabular}
\end{tiny}
\caption{
$(\theta_0,\theta_1,\theta_2,a_1,a_2,a_3,\pi_1,\pi_2,\pi_3)=(6/7,-2/7,3/7,4+i,-1+3i,-2-i,
0.6,0.25,0.15)$. The starting point of left array is
$(\sigma_{ini},(\theta_{ini,i})_{0\leq i \leq 4})=(
0.005,0.2,0.3,0.2,0.2,0.2)$ and the starting point of right array is
$(\sigma_{ini},(\theta_{ini,i})_{0\leq i \leq 4} )=( 0.7, 0.2, -0.3,
0.2,0.002,0.001)$.} \label{tableAR2b}
\end{center}
\end{table}
\bigskip
\pagebreak

$\; $ \newpage

\newpage

\newpage

\newpage

\strut \thispagestyle{empty} \vfill

\pagebreak


\section{Discussion}\label{Conclusion:Sec}
\noindent  {\bf $\bullet$ Interest of our estimation procedure.} Our
procedure estimation does not require a priori more specifications
of the model  than  equation (\ref{model}) then it takes the
advantage to adapt to any situations.

\noindent {\bf $\bullet$ Choice of the Hankel matrix.}
 It would be possible to deal with the Toeplitz matrix instead of the Hankel matrix of the
 $(Z_t)$ since the characterizations of $\theta$ and $\sigma_0$ given by relations (\ref{caractthetasig}) and (\ref{CondBORD}) also hold for the Toeplitz matrix. More generally for
the same reason, it would be possible to consider any
$(p+1)^2$-vector built on the moment of type $(\BBe
((\phi_1(Z_0(s(\xi)))^k\overline{\phi_1(Z_0(s(\xi)))^j}])_{k,j \in
\{1,\ldots,p\}}$ where $\phi_1$ is any complex injective function
defined on $\BBc$. The main difference would lie in the non-trivial
determination of $A^{-1}$, the inverse matrix of $A$ defined in
relation (\ref{EquA}). It also could be extend to some entropy
distance which allows to distinguish variables which have less than
$p$ point of support than the others (Gamboa $\&$ Gassiat [1996]).

\noindent  {\bf $\bullet$ Gaussian noise.} The assumption of a
Gaussian noise is not necessary; actually only an indivisible law is
required for the noise. As previously, it is probably more
complicated to exhibit the matrix $A^{-1}$ which contains among
others the calculations of the $\gamma_{j,k}$.

\noindent  {\bf $\bullet$ Method of  moments.}
Traditionally the method of the moments is not very well appreciated
in estimation problems since a little error on the observations entails a big
error in the final estimation. In our case, since $\sigma_0$ is not
a priori estimated, our method takes intrinsically into account such
a type of error and does not possess the disadvantage of the method
of moments. Nevertheless, the restitution of the distribution of
$(X_t)_{t \in \BBz}$ is based on this method and it would be probably more
efficient to apply the MCMC methods  (see for example Sylvia [2001])
since the application of the inverse filter provides a finite
mixture with a known covariance structure.

\noindent  {\bf $\bullet$ Computational comparison.} In order to
compare empirically our method with already existing numerical
results (Gassiat \& Gautherat [1998], p. 1947), we consider the
second order autoregressive model with a real signal, for $p=2$ and
$\theta=(\theta_0,\theta_1,\theta_2)=(\frac 67,-\frac 27,\frac
37)=(0.8571, -0.2857, 0.4286)$.
\\
All quantities which appear in our estimation procedures, are
adapted  to the real case (in particular the matrix $A$). We
characterize the results obtained in Gassiat \textit{et al.} [1998]
by ${\cal M}_1$ and our  results by ${\cal M}_2$. The algorithm we
implement here consists in selecting randomly a starting point for
which the criteria function $J_n$ takes a value close to zero, and
then the pre-defined function fsolve in Matlab is used to find the
zero. The random selection is made from points $(\theta_0,\theta_1)$
which are uniformly distributed  on the sphere $S_1$. The
signal-to-noise ratio expressed in DB is denoted  SNR.

\begin{center}
\begin{tiny}
\begin{tabular}{|llrc|ccc|}
\hline $\sigma_0$ & SNR& $n$& m\'ethode&
$\widehat{E}(\widehat{\theta}_0) \pm std(\widehat{\theta}_0)$&
$\widehat{E}(\widehat{\theta}_1) \pm std(\widehat{\theta}_1)$
&$\widehat{E}(\widehat{\sigma}_0) \pm std(\widehat{\sigma}_0)$
\\
\hline 0.1 &46&100& ${\cal M}_1$ & $ 0.7500 \pm 0.0647$ &
$-0.1767\pm 0.0668$&$0.0000 \pm 0.0000$ \\
0.1 &46&100& ${\cal M}_2$ & $0.8553 \pm 0.0095$ &$-0.2879 \pm
0.0205$ & $ 0.1004\pm 0.0100$\\
$0.1$& 46&500& ${\cal M}_1$ & $0.7591 \pm 0.0140 $&$-0.1860 \pm
0.0154$ & $ 0.0000\pm 0.0000$ \\
0.1& 46&500& ${\cal M}_2$& $0.8580 \pm 0.0049$&  $-0.2941 \pm
0.0187$ & $0.1035 \pm 0.0031$\\
0.1&46&1000& ${\cal M}_1$&$0.7639 \pm 0.0128$&$-0.1923\pm 0.0142
$&$0.0000 \pm 0.0000$\\
0.1&46&1000& ${\cal M}_2$& $0.8573 \pm 0.0080$& $ -0.2837 \pm
0.0151$ & $0.1009 \pm 0.0139$ \\
\hline \hline
 1 &0.46&500&$ {\cal M}_1$&$ 0.5913\pm 0.2537$ &$-0.0916 \pm
0.2652$ & $0.9407 \pm 0.1015$ \\
1 &0.46&500&$ {\cal M}_2$&$ 0.1423 \pm 0.6017$ & $ -0.4147 \pm
0.2475$&$ 0.8607 \pm 0.5653$ \\
1&0.46&5000&$ {\cal M}_1$ & $ 0.7529\pm 0.1211$ &$-0.2450 \pm  0.1353$ &$ 1.0022 \pm 0.0612 $ \\
 1&0.46&5000&$ {\cal M}_2$ & $ 0.5491 \pm 0.5060$ &$-0.3259 \pm
0.2020$ & $ 0.9646  \pm 0.3915$
\\
1&0.46&15000&${\cal M}_1$ &$0.7771\pm 0.1294$ &$-0.2576 \pm 0.1206$&
$1.0302 \pm  0.0651$ \\
1&0.46&15000&${\cal M}_2$ & $ 0.7159 \pm 0.3757$ & $ -0.2928 \pm
0.1476$ & $ 0.9844 \pm 0.0387$ \\
\hline
\end{tabular}
\end{tiny}
\end{center}

One must note that for a large SNR, the method ${\cal M}_2$ is more
performant than the method ${\cal M}_1$ both in the estimated values
and in the empirical standard deviation. The non-asymptotic side
of our method is particularly highlighted  in the case of a large SNR:
when SNR is equal to $46$, our method  works  well
even for a small sample size and it is always better than the
method ${\cal M}_1$. For a very small SNR
(equal  to $0.46$), the situation is changing:
no method is able to handle such a level of noise;
to our knowledge, only the paper of
Gassiat \& Gautherat [1998] gave numerical results in this situation.
One can note that, even a priori the results do not give satisfactory, our
estimation methods works better and better with  an increasing $n$ . \\
Another aspect which is of importance, is that our method does not
need to calibrate some parameters, on the contrary of the estimation
method in Gassiat \& Gautherat [1998], which is based on the
minimization of a penalized contrast function. Moreover our method
does not require a starting point that we have to fixed in advance
since, as it is mentioned above, the starting point is selected
randomly: it is not the case for the results of the section
\ref{Simu:Sec} which are sensitive to the starting point we choose
nor for the numerical results of the method ${\cal M}_1$ described
in Gassiat \& Gautherat [1998]. The method ${\cal M}_1$ needs more:
the starting point must  be near enough to  the true valeur,
otherwise the method could provide a local minimum. The gain of our
method is that no a priori on extra parameters and on a starting
point is requested, but its drawback  is that in some cases, it
generates a large standard deviation. The outlook of the future use
of the inverse matrix $A^{-1}$ since it is explicit done, would
probably perform the numerical results.

\section{Proofs}\label{Proof:Sec}
For convenience, denote $J(\sigma, s(\xi))=J(\sigma, \xi)$. We first
give a very useful tool which is a combination of existing results
obtained by Gautherat [1997] (see Lemma 5.3.1 p.130) and Gassiat \&
Gautherat [1999] (see Lemma 4.1. p. 1695):
\begin{Lemma}\label{diffHn2} Under assumptions {\bf (M1)}-{\bf (M8)} and {\bf
(P)}, one gets

\begin{description}
\item {\bf i)} $\forall \sigma \in \BBr_+^*, \forall \xi \in
{\cal K},\;\;\; J_n(\sigma,\xi) \xrightarrow[ n \rightarrow \infty
]{\BBp-a.s.} J (\sigma, \xi)$.
\item {\bf ii)} $\forall n \in \BBn, \forall \xi \in {\cal K}$, $ J_n(\sigma,\xi)$ is differentiable with respect to
$\sigma$ and  $\partial_1^1 J_n(\sigma_0,\xi_0) \xrightarrow[ n
\rightarrow \infty ]{\BBp-a.s.}   \partial_1^1 J
(\sigma_0,\xi_0)=-\alpha<0$. \label{diffH1}
\item{\bf iii)}  The function
$s(\xi) \in \Theta \mapsto \partial_1 J_n(\sigma_0,{\bar s}(\xi))$
is continuous on $\Theta$.
\item{\bf iv)} $\forall n \in \BBn$, $ J_n(\sigma,\xi)$ is twice differentiable in $(\sigma_0, \xi_0)$ with respect
to both  $\sigma$ and $\xi$. The first and second derivatives of
$J_n(\sigma,\xi)$ in $\sigma_0$ and $\xi_0$ converge $\BBp-a.s.$ to
the first and second derivative of $J(\sigma,\xi)$ in $(\sigma_0,
\xi_0)$.
\item{\bf v)} The asymptotic distribution of $\left({\sqrt n} \
\partial_2^1 J_n (\sigma_0, \xi_0),{\sqrt n} \;  J_n
(\sigma_0, \xi_0)\right)$ is a centered Gaussian vector variance $$
D^1 h(d(\sigma_0,\xi_0)) A^{-1}(\sigma_0 \|s(\xi_0)\|_2) \Gamma_1
(A^{-1}(\sigma_0 \|s(\xi_0)\|_2))' (D^1 h(d(\sigma_0,\xi_0)))'.$$
\end{description}
\end{Lemma}
\noindent {\bf Proof of Lemma \ref{diffHn2}} These results are
proved using the compactness of ${\cal K}$ and in adapting the proof
of Lemma 4.1 in Gassiat \& Gautherat [1999] to the almost surely
convergence for {\bf i)}-{\bf iv)}, and directly from the statement
adapted to the almost-surely convergence of  Gassiat \& Gautherat
[1999], Lemma 4.1 for {\bf ii)}-{\bf v)}. Whereas, {\bf iii)} is
obtained due  to the truncation $k(n)$ of $s(\xi)$ and the
polynomial structure of
$J_n$. \\


%
\noindent
{\bf Proof of Theorem \ref{consistance}.} \\
{\bf $\bullet $ Consistency of $\widehat{\sigma}_{0}$.} Let ${\cal
V}_{\sigma_0}$ be some neighborhood of $\sigma_0$ and choose
$\sigma_1 < \sigma_0$ in ${\cal V}_{\sigma_0}$.  Due to relation
(\ref{CondBORD}), $\sigma_1$ satisfies $J(\sigma_1, \xi)> 0$ for all
$\xi$. On the other hand,  due to assertion {\bf ii)} of Lemma
\ref{diffHn2} it is always possible to consider $\sigma_2> \sigma_0$
in ${\cal V}_{\sigma_0}$ such that $J(\sigma_2, \xi_0)<0$ and such
that for any $\sigma \in ]\sigma_0, \sigma_2[$,
$J(\sigma,\xi_0)<0$.\\ Assertion {\bf i)} in Lemma \ref{diffHn2}
leads to
\begin{eqnarray*}
J_n ( \sigma_1,\xi_0) & \xrightarrow[ n \rightarrow \infty
]{\BBp-a.s.}  &
J(\sigma_1,\xi_0), \; \\
J_n(\sigma_2,  \xi_0) &\xrightarrow[ n \rightarrow \infty
]{\BBp-a.s.}  & J(\sigma_2,\xi_0).
\end{eqnarray*}
Let $0<\epsilon <\inf \{  J(\sigma_1,\xi_0), |J(\sigma_2,\xi_0)|\}$.
Then,  it exists a positive integer $N_0$ such that for all $n \geq
N_0$, $J_n ( \sigma_1, \xi_0)>0$ and $J_n ( \sigma_2, \xi_0)< 0$.
Thus, from {\bf ii)} in Lemma \ref{diffHn2}, it follows that for all
$n
>N_0$, it exists $\tilde \sigma_n \in ]\sigma_1,
\sigma_2[$ such that $J_n(\tilde \sigma_n,\xi_0)=0$ and
 we choose $\tilde \sigma_n$ such that $\tilde
\sigma_n = \inf \{\sigma_n \in ]\sigma_1, \sigma_2[ \; :\; J_n(
\sigma_n,\xi_0)=0 \}$. From Assertions {\bf i)} and {\bf ii)} in
Lemma \ref{diffHn2}, and a Taylor expansion of $J_n$ at $(\tilde
\sigma_n,\xi_0)$, one obtains,
$$J_n(\tilde \sigma_n,\xi_0) =J_n( \sigma_0,\xi_0) + (\tilde
\sigma_n - \sigma_0) \;\partial_1^1 J_n (\sigma_0,\xi_0) (1+
o(1)),$$
and one gets $(\tilde \sigma_n - \sigma_0) \xrightarrow[ n
\rightarrow \infty ]{\BBp-a.s.}  0$. Since for all $\sigma <
\sigma_1, \;\;J_n(\sigma,\xi_0)
>0$,
$\widehat{\sigma}_{0}$ satisfies $ \widehat{\sigma}_{0}>\sigma_1$.
Consider only large $n$ that is $n$ such that $n>N_0$, by definition
of ${\widehat \sigma}_{0}$, one has ${\widehat \sigma}_{0} < \tilde
\sigma_n$. Since we consider only $\sigma$ lying in the compact set
$[\sigma_1,\sigma_2]$, there exists a subsequence ${\tilde
\sigma}_{0,n}$ of ${\widehat \sigma}_{0}$ which converges to
${\tilde \sigma}_0$ and which satisfies ${\tilde \sigma}_0 <
\sigma_0$. Since $J_n({\tilde \sigma}_{0,n},\xi_0)=0$ and due to
$J_n({\tilde \sigma}_{0,n}, \xi_0) \xrightarrow[ n \rightarrow
\infty ]{\BBp-a.s.} J( {\tilde \sigma}_0, \xi_0)$, it follows that
$J( {\tilde \sigma}_0, \xi_0)=0$, which contradicts the definition
of $ \sigma_0$ (see relation (\ref{CondBORD})). This achieves the
proof.
%
%
%

\noindent {\bf $\bullet $ Consistency of $\widehat{\xi}_0$.}
Consider only $\sigma$'s in $[\sigma_1,\sigma_2]$.  Since ${\cal K}$
is a compact set, ${\widehat \xi}_0$ admits a subsequence $({\tilde
\xi}_n)_n$ which converges to ${\tilde \xi}_0$. Assertion {\bf i)}
in Lemma \ref{diffHn2}, the a.s.-convergence of ${\widehat
\sigma}_{0}$ and the continuity of $J_n$ and $s$ lead to $J_n (
{\widehat \sigma}_{0},{\tilde \xi}_n ) \xrightarrow[ n \rightarrow
\infty ]{\BBp-a.s.} J(\sigma_0,{\tilde \xi}_0).$ This implies that
${\tilde \xi}_0$ is equal to $\xi_0$ since $J(\sigma_0,{\tilde
\xi}_0)=0
\Longleftrightarrow  s({\tilde \xi}_0)=\theta$. \\
Suppose now there exists $\overline {\xi}_0$ an accumulation point
which is different from ${\tilde \xi}_0$.  Then it exists another
subsequence $\overline {\xi}_n$ of $\widehat{\xi}_0$ which converges
to $\overline {\xi}_0$. Using the same tricks as previously, one
gets $s(\overline {\xi}_0)=\theta$ which proves the uniqueness of
${\tilde \xi}_0$.
%

\vspace{0.2cm}
\noindent {\bf  Proof of Corollary  \ref{consistanceapi}.} This
proof is explicitly done in Gautherat [2002]  (see proof of Theorem
3.2). It is only based on the consistency of $\widehat{\sigma}_{0}$
and $\widehat{\xi}_0$.

\vspace{0.2cm}
\noindent {\bf Proof of Theorem \ref{propspeed}.} The definition of
$\widehat{\sigma}_{0}$ leads to $J_n(\widehat{\sigma}_{0},
\widehat{\xi}_0)=0$.  It entails that both $\partial_1^1
J_n^2(\widehat{\sigma}_{0}, \widehat{\xi}_0)=0$ and $\partial^1_{2}
J_n^2(\widehat{\sigma}_0, \widehat{\xi}_0)=0_d$. For simplicity's
sake, denote $J_n(\sigma_0, \xi_0)=~J_n$,   $\partial_{i}^{r_i}
J_n(\sigma_0, \xi_0)=\partial_{i}^{r_i}J_n,\;i=1,2$ and
$\partial_{i,j}^{r_i,r_j}
J_n(\sigma_0,\xi_0)=\partial_{i,j}^{r_i,r_j}J_n$, $i,j=1,2$.
Therefore, one can apply the Delta method to  $J_n^2$ at $(\sigma_0,
\xi_0)$,  since
\begin{eqnarray*}
\left( \begin{array}{c} \partial_2^1 J_n^2(\widehat{\sigma}_{0},
\widehat{\xi}_0)
\\
\partial_1^1
J_n^2(\widehat{\sigma}_{0}, \widehat{\xi}_0)
\end{array}\right)=
\left( \begin{array}{c} 2J_n(\widehat{\sigma}_{0},
\widehat{\xi}_0)\partial_2^1 J_n(\widehat{\sigma}_{0},
\widehat{\xi}_0)
\\
2J_n(\widehat{\sigma}_{0}, \widehat{\xi}_0) \partial_1^1
J_n(\widehat{\sigma}_{0}, \widehat{\xi}_0)
\end{array}\right) =\left( \begin{array}{c} 0_d
\\
0
\end{array}\right).
\end{eqnarray*}
Now, the expansion  at the first order of $J_n^2$ at
$(\sigma_0,\xi_0)$ is
\begin{eqnarray*}
\left(
\begin{array}{cc}
2(\partial_{2}^1J_n)'\partial_{2}^1J_n+2 J_n
\partial_{2}^2J_n
& 2(\partial_{2}^1J_n)' \partial_{1}J_n+2 J_n \partial_{1,2}^{1,1}
J_n
\\
~&~
\\
2\partial_{2}^1J_n \partial_{1}^1J_n+2J_n (\partial_{1,2}^{1,1}
J_n)'&
 2(\partial_{1}^1J_n)^2 +2J_n D^1J_n
\end{array}\right)\!\! \left( \begin{array}{c} \widehat{\xi}_0- {\xi_0} \\
\widehat{\sigma}_{0}-\sigma_0\end{array}\right)\!\!\!(1+o(1)) +\\
\left( \!\!\!\begin{array}{c} 2J_n \partial^1_{2}J_n
\\
2J_n \partial_{1}^1J_n\end{array}\right) \! =\!\!\!\left(
\begin{array}{c} 0_d
\\
0
\end{array}\right)
\end{eqnarray*}
Denote
$A_n=(\partial_{2}^1 J_n)' \partial_{2}^1 J_n+ J_n
\partial_{2}^2 J_n$,
$B_n=(\partial_{2}^1 J_n)' \partial_{1}^1 J_n+J_n
\partial_{1,2}^{1,1} J_n$ and $d_n=(\partial_{1}^1 J_n)^2 +J_n
\partial_1^2 J_n$.
Then, from the Schur complement (Searle [1982]),  we obtain
\begin{eqnarray*}
&\left(\begin{array}{c} \widehat{\xi}_0- {\xi_0}
\\\widehat{\sigma}_0-\sigma_0\end{array}\right)(1+o(1))&\\
& =& \\
&\left( \begin{array}{c}
 J_n A_n^{-1}(\partial_{2}^1J_n)'+ \frac{ J_n}{d_n-B_n'A_n^{-1} B_n}\left(
A_n^{-1}B_nB_n'A_n^{-1} (\partial_{2}^1J_n)'- A_n^{-1}B_n
\partial_{1}^1J_n\right)
\\
\frac{J_n}{d_n-B_n'A_n^{-1} B_n}\left( -B_n'A_n^{-1}
(\partial_{2}^1J_n)'+ \partial_{1}^1J_n\right)
\end{array}\right).& \nonumber
\end{eqnarray*}
For $n$ fixed, note that $J_n(\sigma_0,\xi_0)$ differs from zero.
Thus, one could rewrite  the  up-right term in the previous equation
using the expression of $B_n$ and  dividing it by $J_n$; it gives
three terms $T_1$, $T_2$ and $T_3$ which are defined by:
\begin{eqnarray*}
T_1 &=&
\left(\frac{A_n}{J_n}\right)^{-1}(1-\frac{(\partial_{1}^1J_n)^2}{d_n-B_n'A_n^{-1}
B_n})(\partial_{2}^1J_n)',\\
T_2 & =&
\left(\frac{A_n}{J_n}\right)^{-1}\frac{B_nB_n'A_n^{-1}(\partial_{2}^1J_n)'}{d_n-B_n'A_n^{-1}
B_n},\\
T_3 & = & - \left(\frac{A_n}{J_n}\right)^{-1}
\frac{\partial_{1}^1J_n}{d_n-B_n'A_n^{-1} B_n}J_n
\partial_{1,2}^{1,1}J_n.
\end{eqnarray*}
Rewrite the approximation of the vector $(\widehat{\xi}_0- {\xi_0},
\widehat{\sigma}_{0}-\sigma_0)'$ as follows
\begin{eqnarray*}
J_n   \! \! \left(
\begin{array}{c}
  \! \! \!  T_1+T_2+T_3
\\
\frac{1}{d_n-B_n'A_n^{-1} B_n}\left( -B_n'A_n^{-1}
(\partial_{2}^1J_n)'+ \partial_{1}^1J_n\right)
\end{array}\right).
\end{eqnarray*}
Due to Lemma \ref{diffH1} and the continuity of $J_n$ in $\xi$,
one has $A_n \xrightarrow[n \rightarrow \infty]{\BBp-a.s.} 0$,
$B_n \xrightarrow[n \rightarrow \infty]{\BBp-a.s.}  0$, $d_n
\xrightarrow[n \rightarrow \infty]{\BBp-a.s.} (\partial_{1}^1
J(\sigma_0,\xi_0))^2=\alpha^2>0$, $\frac{A_n}{J_n} \xrightarrow[n
\rightarrow \infty]{\BBp}
\partial_{2}^2 J(\sigma_0,\xi_0)$,
$\frac{B_n}{J_n}  \xrightarrow[n \rightarrow \infty]{{\cal L}} W$,
where $W$ is a $d$-dimensional non degenerate random vector and
$\frac{d_n}{J_n} \xrightarrow[n \rightarrow \infty]{\BBp-a.s.}
\pm\infty$.
As $n$ large enough,  it entails that
\begin{eqnarray}
\sqrt{n} \left( \begin{array}{c}  \widehat{\xi}_0- {\xi_0}
\\\widehat{\sigma}_{0}-\sigma_0\end{array}\right)
& \left. \begin{array}{c} {{\cal L}} \\
 =
\end{array} \right.&  \sqrt{n} J_n(\sigma_0,\xi_0) \left(
\begin{array}{c}
\frac{(\partial_{2}^2J(\sigma_0,\xi_0))^{-1}(\partial_{1,2}^{1,1}
J(\sigma_0,\xi_0))}{\alpha}
\\
~ \\
\frac{1}{\alpha}
\end{array}\right). \label{TCL}
\end{eqnarray}
Note that $h({\tilde d}(\sigma_0, \xi_0))= J_n(\sigma_0,\xi_0),$
where $h$ is the determinant function. Then, due to Assumption {\bf
(M7)} and due to the Taylor expansion of $h$ at ${\tilde
d}(\sigma_0, \xi_0)$, one obtains
\begin{eqnarray}
h({\tilde d}(\sigma_0, \xi_0))& =&D^1h(\tilde d(\sigma_0,\xi_0))
A^{-1}(\sigma_0\|\theta\|_2) (d_n(\xi_0)-d(\xi_0)+
o(\frac{1}{\sqrt{n}})). \label{Taylor}
\end{eqnarray}
Set $M= A^{-1}(\sigma_0 \|s(\xi_0)\|_2) \Gamma_1 (A^{-1}(\sigma_0
\|s(\xi_0)\|_2))' $ and \\
$N=\left(
\begin{array}{c}
\frac{1}{\alpha} \; (\partial_{2}^2J (\sigma_0,\xi_0))^{-1} \;
\partial_{1,2}^{1,1}J (\sigma_0,\xi_0)
\\
~ \\ \frac{1}{\alpha}
\end{array}\right)$, then due to (\ref{TCL}), (\ref{Taylor}) and {\bf v)} of Lemma  \ref{diffH1} one gets,
\begin{eqnarray*}
\sqrt{n}\left( \begin{array}{c} \widehat{\xi}_0- \xi_0
\\
\widehat{\sigma}_{0}-\sigma_0\end{array}\right)\xrightarrow[n
\rightarrow \infty]{{\cal L}} {\cal N} \left(0_{d+1}, N
D^1h(\tilde d(\sigma_0,\xi_0)) \; M \;(D^1h(\tilde
d(\sigma_0,\xi_0)))'N'\right) \label{Distrib:thetasigma}.
\end{eqnarray*}

\vspace{0.2cm}
\noindent {\bf  Proof of Corollary \ref{propspeedapi}.} Following
the proof of Theorem 3.3 in Gautherat [2002], it remains to obtain
an equivalent for $\sqrt{n} ({\tilde d}_n({\hat
\sigma}_{0},{\widehat \xi}_0)-{\tilde d}(\sigma_0,\xi_0))$. As $n$
large enough, this term is equivalent in distribution to
%
%
 $$\left( {\rm Id}_{(p+1)^2} + D^1 {\tilde d}(\sigma_0,\xi_0) D^1
 h(d(\sigma_0,\xi_0)) N \right) A^{-1}(\sigma_0,\|\theta\|_2) \sqrt{n}
(d_n(\xi_0)-d(\xi_0)).$$
%
%
On the other hand, one has
\begin{eqnarray}
\sqrt{n} ({\widehat a}-a) &=& \left(\frac{C^{-1}}{2
|v_p^*|^2}B\right) \sqrt{n} ({\tilde d}_n({\hat \sigma}_{0},\widehat
\xi_0)-{\tilde
d}(\sigma_0,\xi_0)),  \label{TCLa} \\
\sqrt{n} ({\widehat \Pi}-\Pi)&=& L^{-1}\left( Proj + F
\frac{C^{-1}}{2 |v_p^*|^2}B\right) \sqrt{n} ({\tilde d}_n({\widehat
\sigma}_{0},\widehat \xi_0)-{\tilde d}(\sigma_0,\xi_0)).
\label{TCLpi}
\end{eqnarray}
Relations (\ref{TCLa}) and (\ref{TCLpi}) entail the results. All
matrices used here are defined in the statement of both  Theorem
\ref{propspeed} or  Corollary \ref{propspeedapi}.
\bibliographystyle{acm}
\bibliography{bibbase}
%
%
%
%
%


\end{document}